\documentclass[preprint,11pt,3p]{elsarticle}

\usepackage{amssymb}
 \usepackage{amsthm}
\usepackage{color}
\usepackage{subfig}
\usepackage{graphicx}
\usepackage{amsmath,amsfonts}
\usepackage{bm}
\usepackage{soul}
\usepackage{sidecap}
\usepackage{ifthen}
\usepackage{comment}

  \newtheorem{thm}{Theorem}
 \newtheorem{lem}[thm]{Lemma}
 \newdefinition{defn}{Definition}
 \newproof{pf}{Proof}

\journal{\empty}

\begin{document}

\begin{frontmatter}



\title{Landau-Lifshitz Equation with Affine Control}

 \cortext[cor1]{Corresponding author}
\author[UW]{Amenda N. Chow\corref{cor1}}
\ead{a29chow@uwaterloo.ca}
 \author[UW]{Kirsten A. Morris}
 \ead{kmorris@uwaterloo.ca}
\address[UW]{Department of Applied Mathematics, University of Waterloo, Canada}

\begin{abstract}
The Landau-Lifshitz equation is a coupled set of nonlinear partial differential equations that describes the dynamics of magnetization in a ferromagnet. This equation has an infinite number of stable equilibria.  Steering the system from one  equilibrium to another is a problem of both theoretical and practical interest.  Since the objective is to steer between equilibria, approaches based on linearization are not appropriate.  It is proven that affine  proportional control can be used to steer the system from an arbitrary initial state, including an  equilibrium point, to a specified equilibrium point.  The second point becomes a globally asymptotically stable  equilibrium of  the controlled system. The control also removes  hysteresis from the Landau-Lifshitz equation.  These results are illustrated with simulations.
\end{abstract}


\end{frontmatter}


\section{Introduction}
The Landau-Lifshitz equation was  developed to model the behaviour of domain walls in magnetic regions within ferromagnetic structures \cite{Landau1935}. For example, the one-dimensional Landau-Lifshitz equation can be used to describe ferromagnetic nanowires, which are often found in memory storage devices such as hard disks, credit cards or tape recordings. Each set of data stored in a memory device is uniquely assigned to a specific stable magnetic state of the ferromagnet.  This can be difficult to achieve due to the presence of hysteresis. Hysteresis is characterized by the presence of multiple equilibria, and  looping in the input-output map is typical \cite{Chow2013ACC,Morris2011}. Consequently, a particular input can lead to different magnetizations.  Therefore, it is desirable to  control magnetization between different stable equilibria.

The Landau-Lifshitz equation is known to exhibit hysteretic behaviour.  For example, \cite{Cowburn1999,Suess2002}  investigated via experiments the shape change  of the hysteresis loop as the structure of the nanomagnet is varied. Experiments conducted on nanowires also demonstrate hysteresis loops \cite{Noh2012}.  Numerical simulations illustrating hysteresis loops are found in  \cite{Wiele2006,Yang2011}. The dynamics of hysteresis  in the Landau-Lifshitz equation has also been represented by a hysteresis operator \cite{CarbouEfendiev2009,Visintin1997}. In much of the aforementioned literature, the presence of hysteresis in the Landau-Lifshitz equation is identified by the fact that input--output curves exhibit a looping behaviour.  This alone is not enough to characterize hysteresis \cite{Chow2013ACC,Morris2011,Bernstein2005}. A looping behaviour must persist with low frequency periodic inputs.
\begin{defn}\label{defBernstein}{\em \cite{Bernstein2005}}
A system exhibits {\it hysteresis} if  a nontrivial closed curve in the input--output map persists for a periodic input as the frequency component of the input signal approaches zero. 
\end{defn}
Another approach to hysteresis is based on the existence of multiple stable equilibria, which are present in the (uncontrolled) Landau-Lifshitz equation \cite[Chapter~6]{Guo2008}.
 \begin{defn}{\em \cite[Definition~3]{Morris2011}} \label{defmultiequilibrium}
A system exhibits {\it hysteresis} if it has \\
(a) multiple stable equilibrium points and\\
(b) dynamics that are considerably faster than the time scale at which inputs are varied.
\end{defn}
 Note that condition (b) is relative to the speed at which a controlled input  is changed.  In many cases, hysteresis is present but is rate-dependent \cite{Morris2011}. 

There is now an extensive body of results on control and stabilization of linear partial differential equations  (PDEs); see for instance the books \cite{Bensoussan-book,Curtain1995,LT00_1,LT00_2} and the review paper \cite{Morris_control_handbook_rev2010}. There are fewer results on control and stabilization of nonlinear partial differential equations  and the Landau-Lifshitz equation is particularly problematic.  Stability of the Landau-Lifshitz equation is often based on linearization  \cite{CarbouLabbe2006,Carbou2006,Carbou2011,Jizzini2011,Labbe2012}.  Local asymptotic stability is shown in \cite{Chow2015} for the controlled linearized Landau-Lifshitz equation. However, because the Landau-Lifshitz equation is not quasi-linear,  analysis based on a linearization may not predict stability of the original system; see \cite{CoronNguyen2015,alJamal2013,AM2014}. Also, when the goal is to steer between equilibria, global stability results are needed. Experiments and numerical simulations on the control of domain walls in a nanowire are presented in \cite{Noh2012,Wieser2011}.  In  \cite{Carbou2008,Carbou2009}, solutions to the Landau-Lifshitz equation are shown to be arbitrarily close to domain walls given a constant control. 

In the next section, the uncontrolled Landau-Lifshitz equation and its equilibrium points are described.  In \cite{Chow2013ACC}, simulations were used to show the Landau-Lifshitz and the linearized Landau-Lifshitz equation exhibit hysteresis. This suggests hysteresis is not due entirely to nonlinearity. These results are reviewed in Section~\ref{secLL}.  Theorem~\ref{thmzeroeigenvalue} demonstrates the linearized uncontolled Landau-Lifshitz equation has a zero eigenvalue.  This suggests use of a proportioonal  controller to stabilize the equation about a given point. It is then proven in  Section~\ref{seccontrol} that   stabilization of the full Landau-Lifshitz equation is achieved with a proportional affine control.  Proportional control can be used to steer the system to an arbitrary equilibrium point of the uncontrolled equation;  in fact, the system can be steered between these points. Simulations illustrating these results are presented in Section~\ref{secExample}.    
  Furthermore, simulations  indicate that hysteresis is absent in the controlled system.   
  
\section{Landau-Lifshitz Equation and Hysteresis}\label{secLL}
Consider the magnetization 
$$
\mathbf m(x,t)=(m_1(x,t),m_2(x,t),m_3(x,t)),
$$  
at position $x \in [0,L]$  and time $t \geq 0 $
in a long thin ferromagnetic material of length $L>0$.  If only the exchange energy term is considered,
 the magnetization
 is modelled by the one--dimensional (uncontrolled) Landau-Lifshitz equation 
 \cite{Brown1963},\cite[Chapter~6]{Guo2008}
 \begin{subequations} \label{eqLLcomplete}
\begin{align}\label{eqLLGuoDing}
\frac{\partial \mathbf m}{\partial t} &= \mathbf m \times  \mathbf m_{xx}-\nu\mathbf m\times\left(\mathbf m\times\mathbf m_{xx}\right)\\
\mathbf m(x,0)&=\mathbf m_0(x)\\
 \mathbf m_x(0,t)&=  \mathbf m_x(L,t)=\mathbf 0.\label{eqboundarycondition}
\end{align}
\end{subequations}
where $\times$ denotes the cross product and $\nu\geq  0$ is the damping parameter, which depends on the type of ferromagnet. The notation $\mathbf m_{x}$  and $\mathbf m_{xx}$  means the magnetization is differentiated with respect to $x$ once and twice, respectively. The gyromagnetic ratio multiplying $ \mathbf m \times  \mathbf m_{xx}$ has been normalized to $1$ for simplicity.  For details on the damping parameter and gryomagnetic ratio, see \cite{Gilbert2004}.  It is assumed there is no magnetic flux at the boundaries and so  Neumann boundary conditions are appropriate.

Define  $\mathcal L_2^3 = \mathcal L_2 ([0,L]; \mathbb R^3)$ with the usual inner product and norm, denoted $\|\cdot\|_{\mathcal L_2^3}$, and  the operator
\begin{equation}
f(\mathbf m)=\mathbf m \times  \mathbf m_{xx}-\nu\mathbf m\times\left(\mathbf m\times \mathbf m_{xx}\right), 
\label{defn:f}
\end{equation}
and its domain 
\begin{align}
\label{setDforfullLL}
D=\{ \mathbf m\in \mathcal L_2^3 :  \mathbf m_x \in \mathcal L_2^3,  \mathbf m_{xx} \in \mathcal L_2^3, 
 \mathbf m_x(0)=\mathbf m_x(L) = \mathbf 0    \}.  
\end{align}
The following theorem is a consequence of the existence and uniqueness results in \cite{Carbou2001,Alouges1992}.
\begin{thm} \label{thmuncontrolledLLsemigroup} If $m(x,0) \in \mathcal L_2^3,$ then the operator $f(\mathbf m)$ with domain  $ D$ generates a nonlinear contraction semigroup on $\mathcal L_2^3.$
\end{thm}

Ferromagnets  are magnetized to saturation \cite[Section~4.1]{Cullity2009}; that is
\[
  || \mathbf m_0 (x) ||_{2} =M_s
\]
where  $||\cdot||_{2}$ is the Euclidean norm and $M_s$ is the magnetization saturation.  Physically,  this  means that at each point, $x$, the magnitude of $\mathbf m_0(x)$ equals the magnetization saturation.  In much of the literature, $M_s$ is set to $1$; see for example, \cite[Section~6.3.1]{Guo2008},  \cite{Carbou2001,Alouges1992,Lakshmanan2011}. This convention is used here. The magnitude of the magnetization does not change with time.
\begin{thm}{\em \cite[Lemma~6.3.1]{Guo2008}}
If $||\mathbf m_0(x)||_{2}=1$, then for all $t\geq0$ the solution to (\ref{eqLLGuoDing}) satisfies 
 \begin{equation}\label{eqconstraint}
  || \mathbf m(x,t)||_{2} =1 .
 \end{equation} 
\end{thm}
\noindent 
The initial condition $\mathbf m_0(x)$ is assumed to be real--valued, which implies $m(x,t)$ is real-valued for all time.

The set of equilibrium points of (\ref{eqLLcomplete}) is  \cite[Theorem~6.1.1]{Guo2008}
\begin{align}\label{equilibriumset}
E=&\{\mathbf a=(a_1,a_2,a_3) : a_1,a_2,a_3 \mbox{ constants and }\mathbf a^\mathrm{T}\mathbf a=1 \}.
\end{align} 
In  \cite[Proposition~6.2.1]{Guo2008}, the stability of the equilibrium points is established using Lyapunov's Theorem and the Lyapunov function 
\[
V(\mathbf m )=\frac{1}{2}\left| \left|  \mathbf m_{x}\right|\right|_{\mathcal L_2^3}^2.
\] Furthermore, $E$ is an asymptotically stable equilibrium set, as stated in the following theorem. 
The proof is the same as that in \cite[Proposition~6.2.1]{Guo2008} except it is for equilibrium sets, rather than equilibrium points. However, individual equilibrium points are only stable, not asymptotically stable. Control is needed to obtain asymptotic stability as illustrated in Section~\ref{seccontrol}. 
\begin{thm}\label{thmlyapunov}
The equilibrium set in (\ref{equilibriumset}) is asymptotically stable in the $\mathcal L_2^3$--norm.
\end{thm}

\begin{figure}\hspace*{-0.3cm}
\centering\scalebox{0.545}{
    \subfloat[{\Large $\omega=1$}]{ \includegraphics{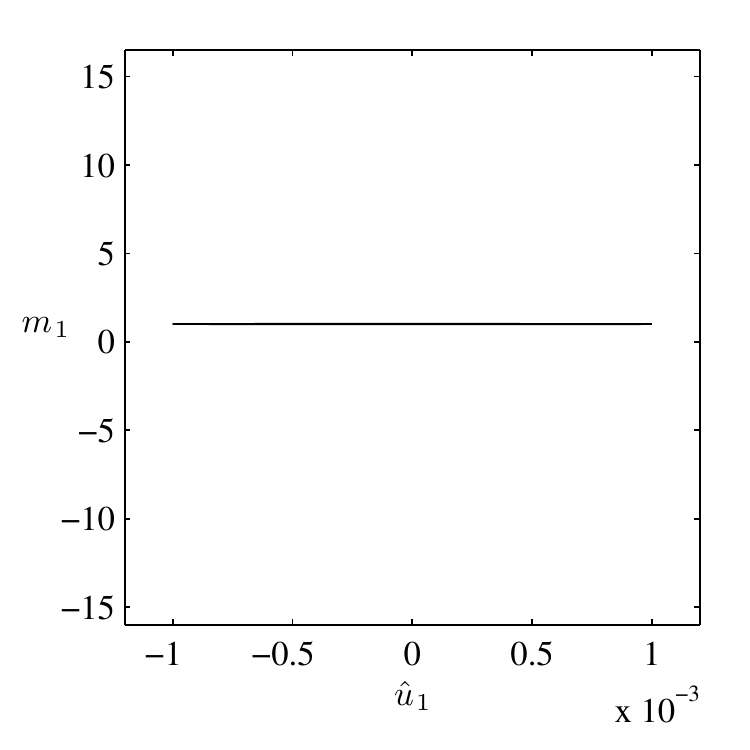}}            
  \subfloat[{\Large$\omega=0.1$}]{\includegraphics{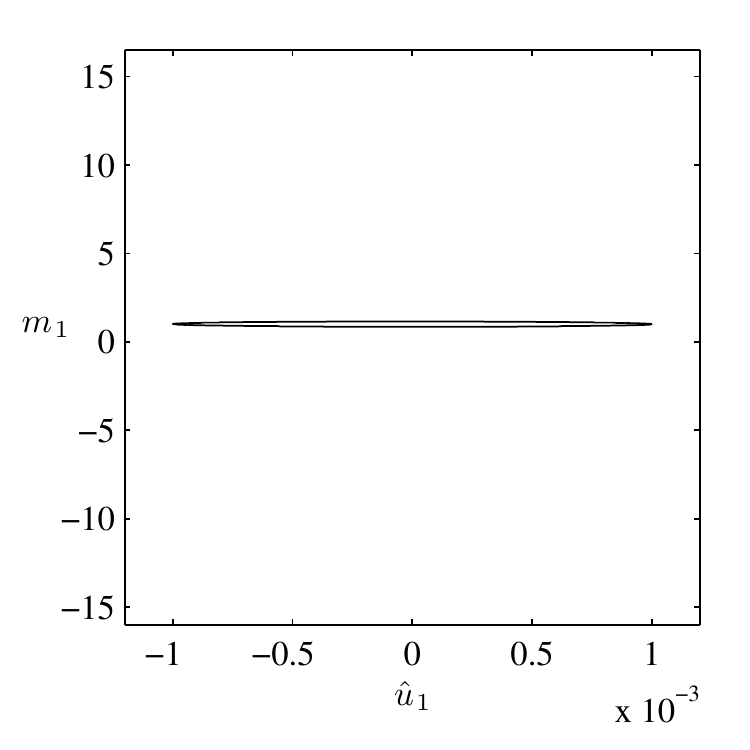}}
 \subfloat[{\Large$\omega=0.01$}]{ \includegraphics{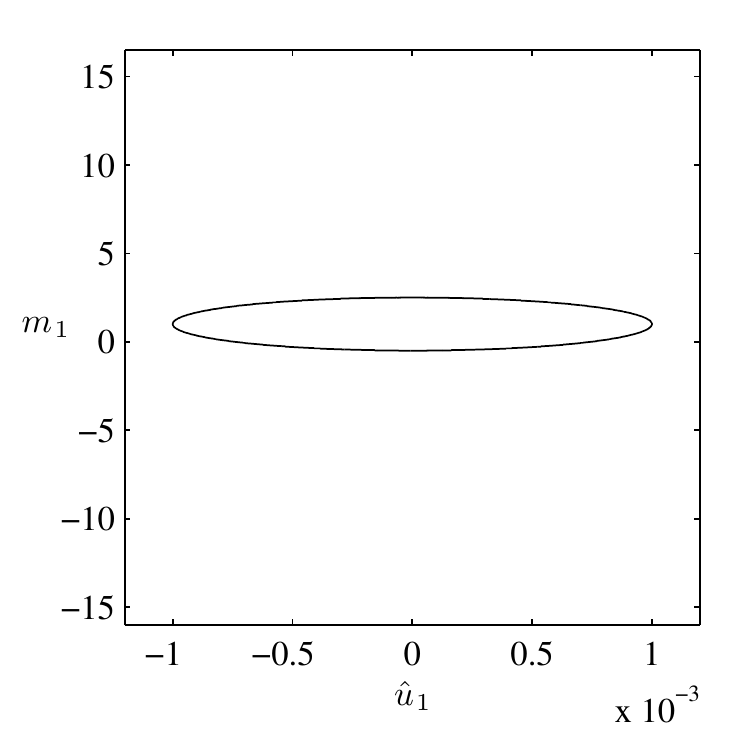}}  
    \subfloat[{\Large$\omega=0.001$}]{ \includegraphics{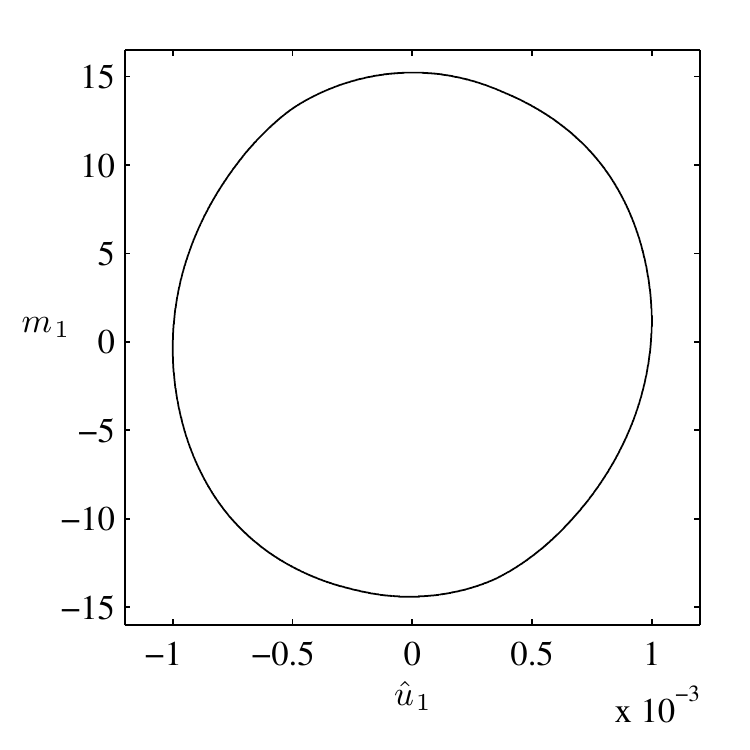}}}
\hspace*{-0.2cm}\scalebox{0.545}{
   \subfloat[{\Large$\omega=1$}]{ \includegraphics{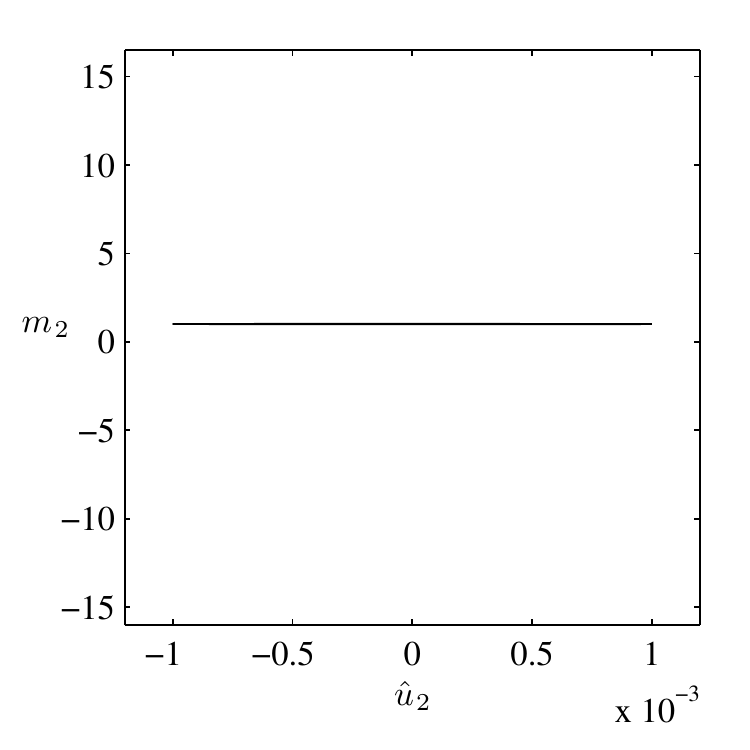}}           
  \subfloat[{\Large$\omega=0.1$}]{\includegraphics{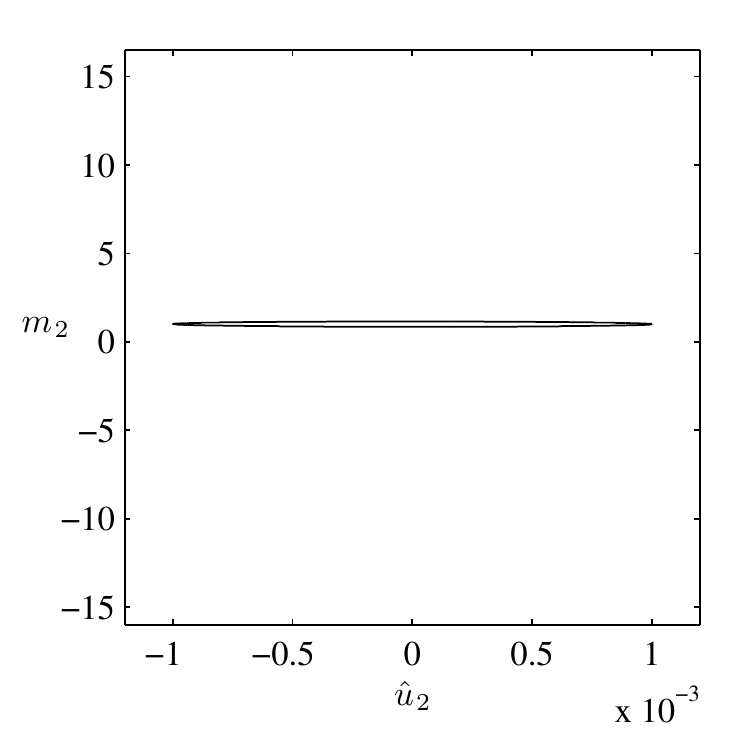}}
   \subfloat[{\Large$\omega=0.01$}]{ \includegraphics{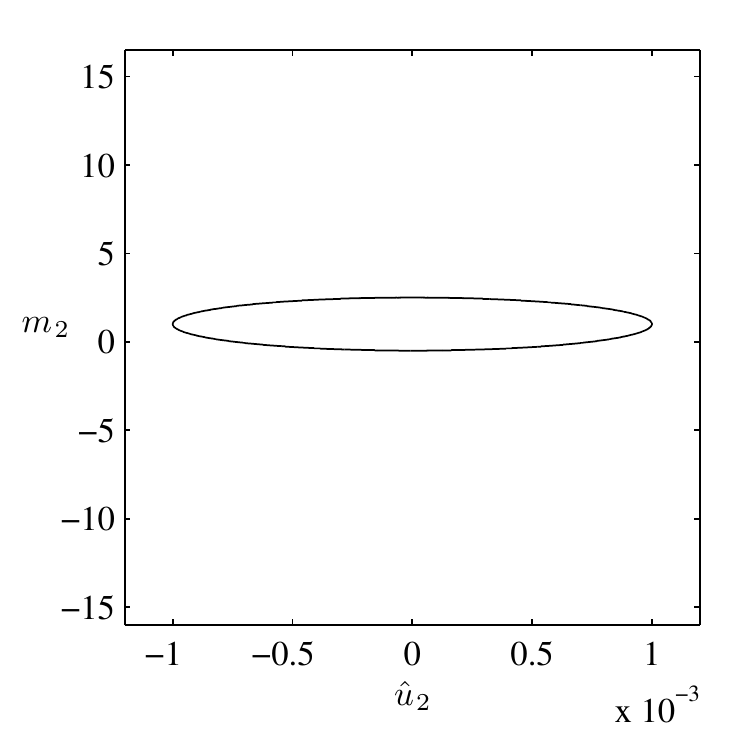}}  
  \subfloat[{\Large$\omega=0.001$}]{ \includegraphics{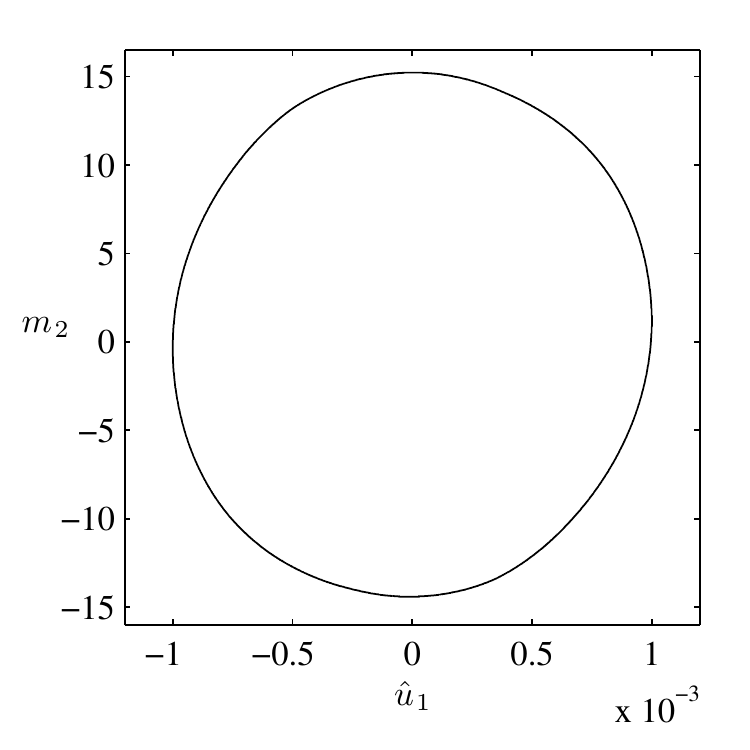}}}
\hspace*{-0.2cm}\scalebox{0.545}{
   \subfloat[{\Large$\omega=1$}]{ \includegraphics{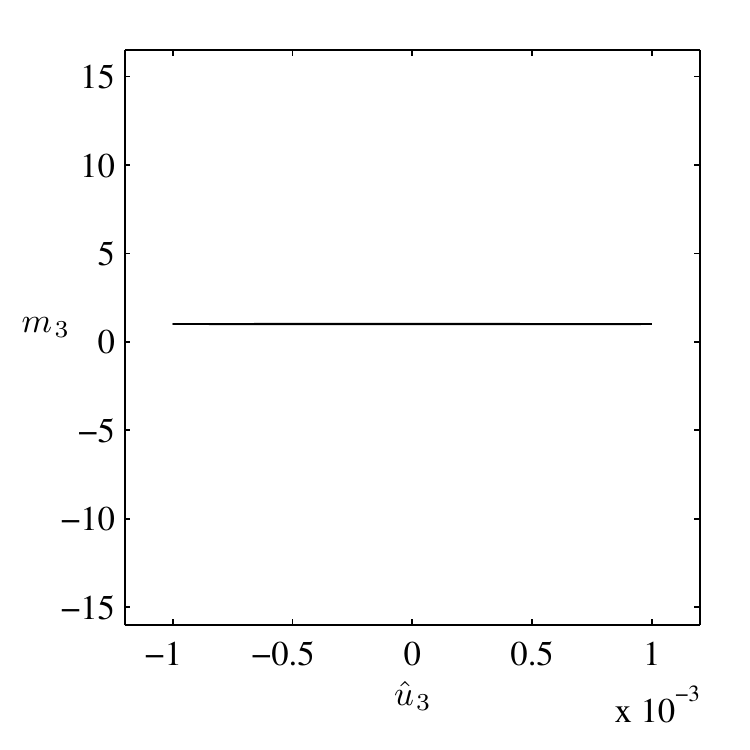}}
   \subfloat[{\Large$\omega=0.1$}]{\includegraphics{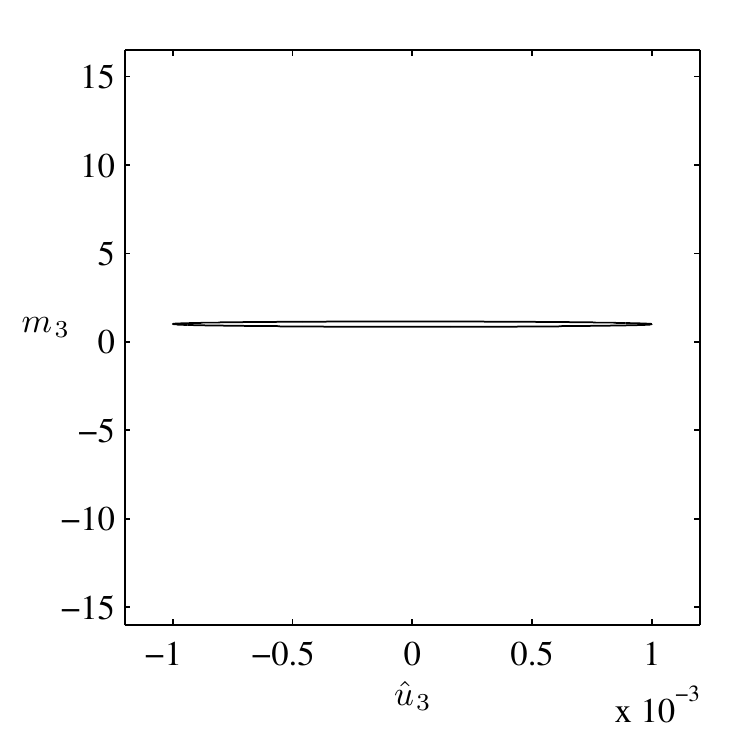}}
   \subfloat[{\Large$\omega=0.01$}]{ \includegraphics{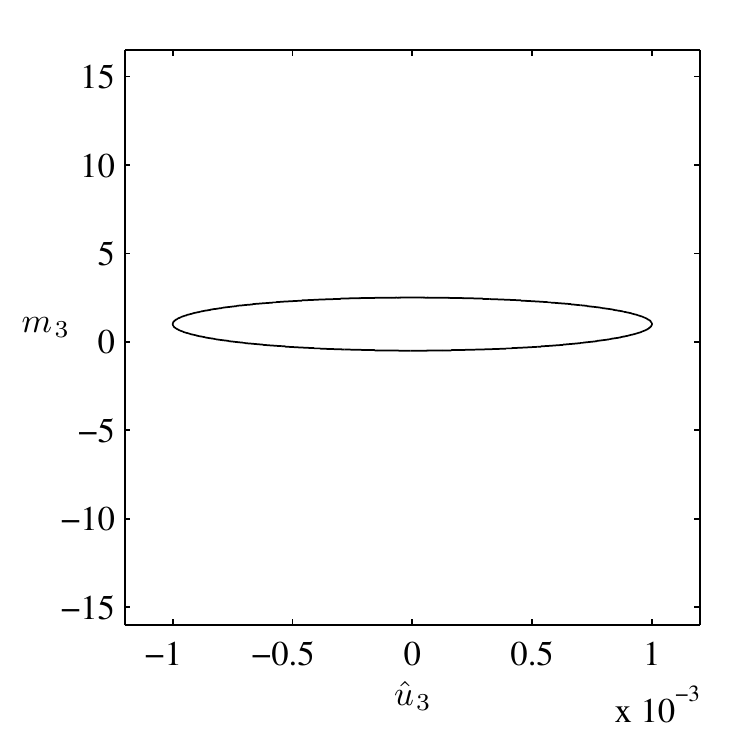}}
  \subfloat[{\Large$\omega=0.001$}]{ \includegraphics{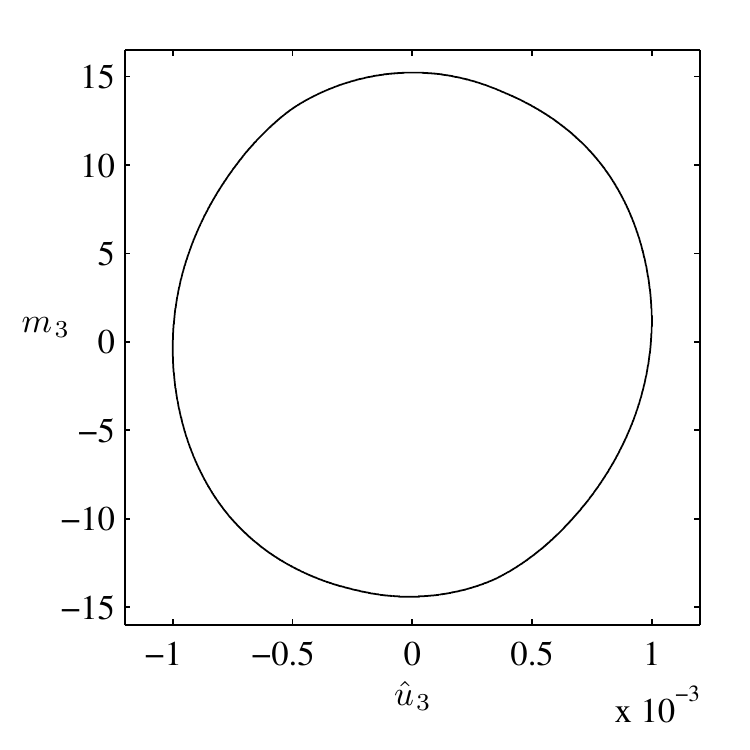}}   }
\caption{\label{figLoopFrequencyVaried} \small  Input--output curves for the (nonlinear) Landau-Lifshitz equation demonstrate persistent looping behaviour as the frequency of the periodic input, $\hat{\mathbf u}$, approaches zero and hence indicates the presence of hysteresis. (a)--(d) Input--output curves for $m_1(x,t)$ with $\hat{\mathbf u}=\left(0.001\cos(\omega t),0,0\right)$ and $\mathbf m_0(x)=\left(1,0,0\right)$. (e)--(h) Input--output curves for $m_2(x,t)$ with  $\hat{\mathbf u}=\left(0,0.001\cos(\omega t),0\right)$ and $\mathbf m_0(x)=\left(0,1,0\right)$. (i)--(l) Input--output curves for $m_3(x,t)$ with $\hat{\mathbf u}=\left(0,0,0.001\cos(\omega t)\right)$ and $\mathbf m_0(x)=\left(0,0,1\right)$. ($L=1$, $\nu=0.02$,  $x=0.6$)} 
\end{figure}

The existence of multiple stable equilibria indicates the presence of hysteresis in the Landau-Lifshitz equation (care of Definition~\ref{defmultiequilibrium}).
 Definition~\ref{defBernstein} is used to establish hysteresis in simulations of the Landau-Lifshitz equation.  For the simulations, a Galerkin approximation for the Landau-Lifshitz equation using  linear spline elements is used. The number of elements is 5 and a periodic input, $\hat{\mathbf u}=\left(0,0.001\cos(\omega t),0\right)$, is applied to the Landau-Lifshitz equation to construct the input-output map.
 Plots of $\mathbf m(x,t)$ with $x$ fixed  against the periodic input are illustrated in Figure~\ref{figLoopFrequencyVaried} for varying frequencies $\omega$.    It is clear from Figure~\ref{figLoopFrequencyVaried} the input--output curves exhibit persistent looping behaviour as the frequency of the input approaches zero.  The continuum of equilibrium points explains the absence of sharp jumps that often appear in hysteresis loops. The similar appearance of the loop shapes between  $m_1(x,t), m_2(x,t), m_3(x,t)$ is due to the symmetric structure of the Landau-Lifshitz equation.

To obtain the linear uncontrolled Landau-Lifshitz equation, equation (\ref{eqLLGuoDing}) is first rewritten in semilinear form,
\begin{equation}\label{eqLLGuoDingTwo}
\frac{\partial \mathbf m}{\partial t} =\nu  \mathbf m_{xx}
 +\mathbf m \times  \mathbf m_{xx}
+\nu\left| \left|  \mathbf m_{x}\right| \right|_{2}^2 \mathbf m,
\end{equation}
using equation~(\ref{eqconstraint}) and properties of cross products, and then $\mathbf m (x,t) = {\mathbf a}+ \mathbf z (x,t)$ is substituted into (\ref{eqLLGuoDingTwo})
 where $\mathbf a \in E$ is an equilibrium of (\ref{eqLLcomplete}) and $\mathbf z \in \mathcal L_2^3$ is a small perturbation. The Landau-Lifshitz equation linearized about an equilibrium $\mathbf a$ is
 \begin{subequations}\label{eqstatespaceform}
\begin{align}\dot{ \mathbf z}&=A\mathbf z,\quad \mathbf z(0)=\mathbf z_0 \\
\mathbf  z_{x}(0,t)&=\mathbf 0=\mathbf  z_{x}(L,t) 
\end{align}
\end{subequations} 
where $A$ is the linear operator, $$A\mathbf z=  \nu\mathbf z_{xx}+\mathbf a \times \mathbf {z}_{xx},$$
and the domain is
\begin{align*}
D(A)= \{ \mathbf z: \mathbf z\in \mathcal L_2^3, \,\,  \mathbf  z_{x} \in \mathcal L_2^3, \, \,  \mathbf z_{xx}^{} \in \mathcal L_2^3, \mathbf  z_{x}(0,t)=\mathbf 0=\mathbf  z_{x}(L,t)\},
\end{align*}
Using \cite[Theorem~6.2]{Banks2012}, the linear operator $A$ can be shown to generate an analytic semigroup; for details, see \cite[Theorem~4.16]{Amenda_thesis}. 

\begin{thm}\label{thmzeroeigenvalue}\cite[Theorem~5]{Chow2013ACC}
Any constant $\mathbf c \in \mathbb R^3$ is a stable equilibrium of (\ref{eqstatespaceform}).
\end{thm}

\begin{pf}
For completeness, the proof is included here. 
Since $A$ generates an analytic semigroup, the spectrum determined growth assumption is satisfied and so the eigenvalues of $A$ determine the stability of the linear system (\ref{eqstatespaceform}) \cite[Section~5.1]{Curtain1995}, \cite[Section~3.2]{Luo1999}.

It is clear that any constant function $\mathbf c$ is an equilibrium of (\ref{eqstatespaceform}). Let $\lambda \in \mathbb C$. The eigenvalue problem of (\ref{eqstatespaceform}) is $\lambda\mathbf v=A\mathbf v$ and  
boundary conditions $\mathbf v_x(0)=\mathbf v_x(L)=\mathbf 0$ where $\mathbf v \in \mathcal L_2^3$. 
Solving, the eigenvalues of (\ref{eqstatespaceform}) are the zero eigenvalue, $\lambda_1=0,$ which is associated to a nonzero constant eigenvector, and the remaining eigenvalues are of the form 
  \begin{align*}
 \lambda_2^{+,-} &= \frac{-(1+2n)^2\pi^2\nu}{L^2}\pm i\frac{(1+2n)^2\pi^2}{L^2}, \qquad
       &  \lambda_3 =  \frac{-(1+2n)^2\pi^2\nu}{L^2}, \\
     \lambda_4^{+,-} &= \frac{-(2n)^2\pi^2\nu}{L^2}\pm i\frac{(2n)^2\pi^2}{L^2}, \qquad
           & \lambda_5 = \frac{-(2n)^2\pi^2\nu}{L^2}
               \end{align*}
where $n\in \mathbb Z$. Since all the eigenvalues have nonpositive real part, the equilibria of (\ref{eqstatespaceform}) are stable.  \qed
\end{pf}
 
\begin{figure}[h]\hspace*{-0.7cm}
\centering\scalebox{1}{
    \subfloat[$\omega=1$]{ \includegraphics[trim = 0mm 60mm 0mm 60mm, clip,width=0.26\textwidth]{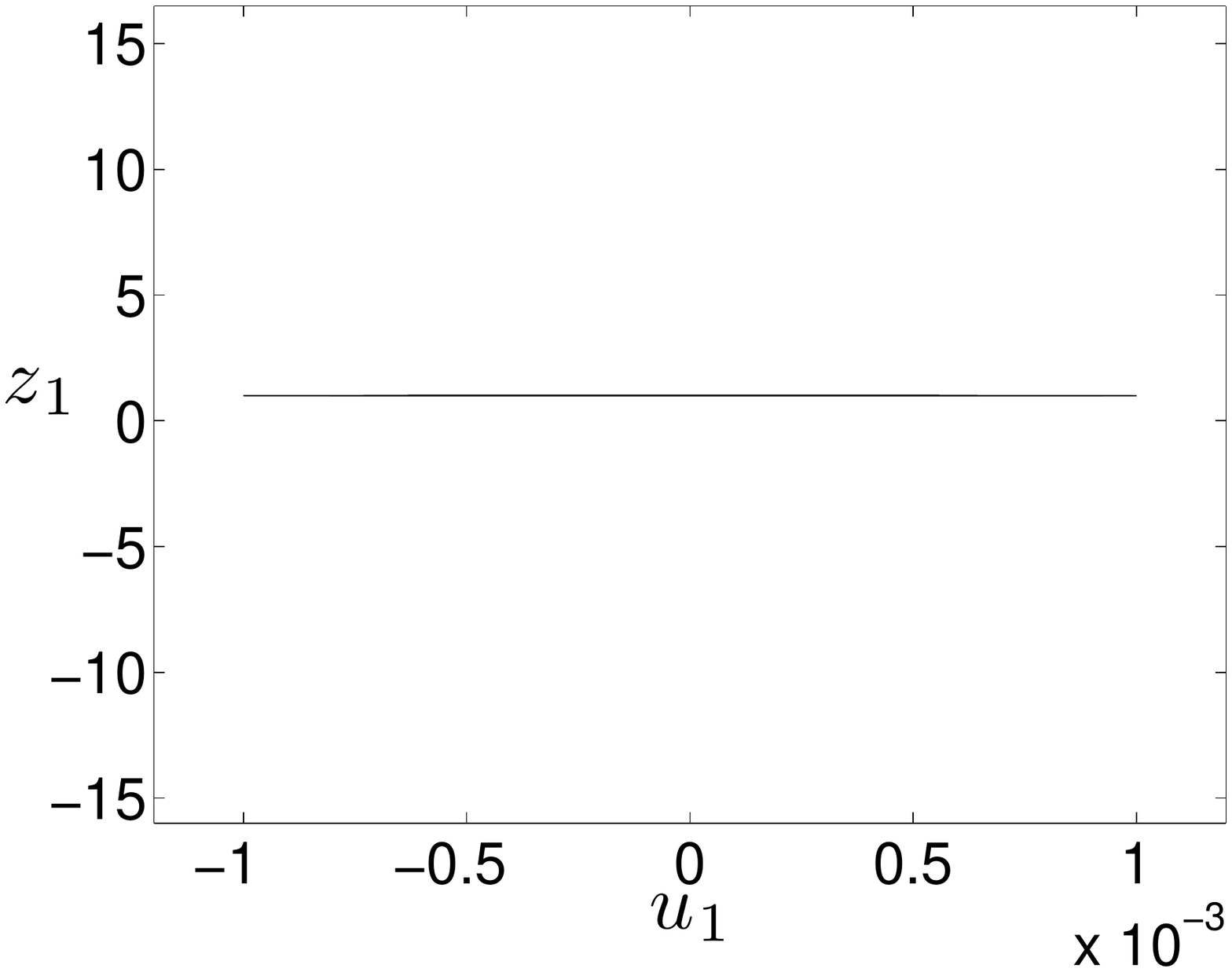}} 
  \subfloat[$\omega=0.1$]{\includegraphics[trim = 0mm 60mm 0mm 60mm, clip,width=0.26\textwidth]{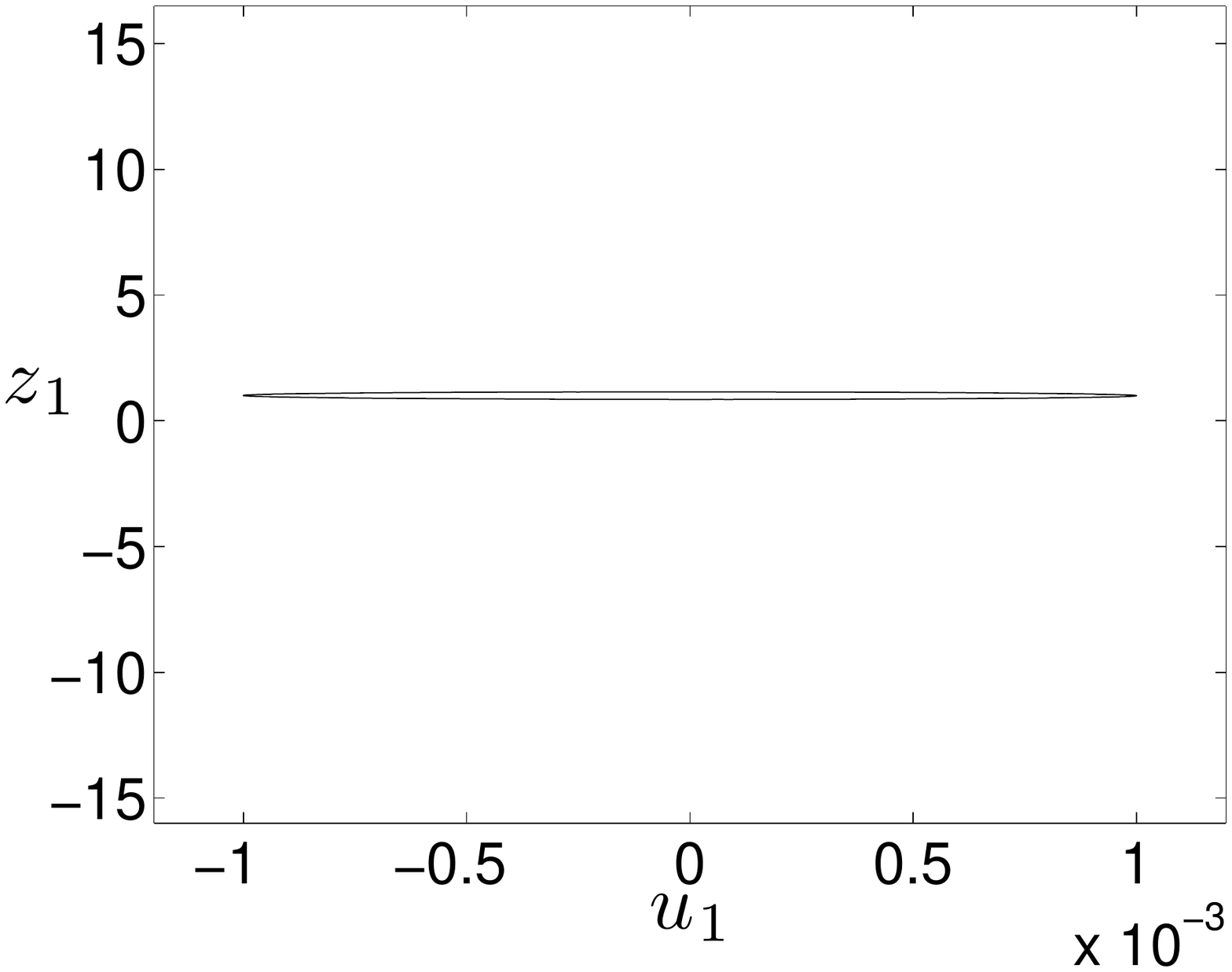}}
  \subfloat[$\omega=0.01$]{ \includegraphics[trim = 0mm 60mm 0mm 60mm, clip,width=0.26\textwidth]{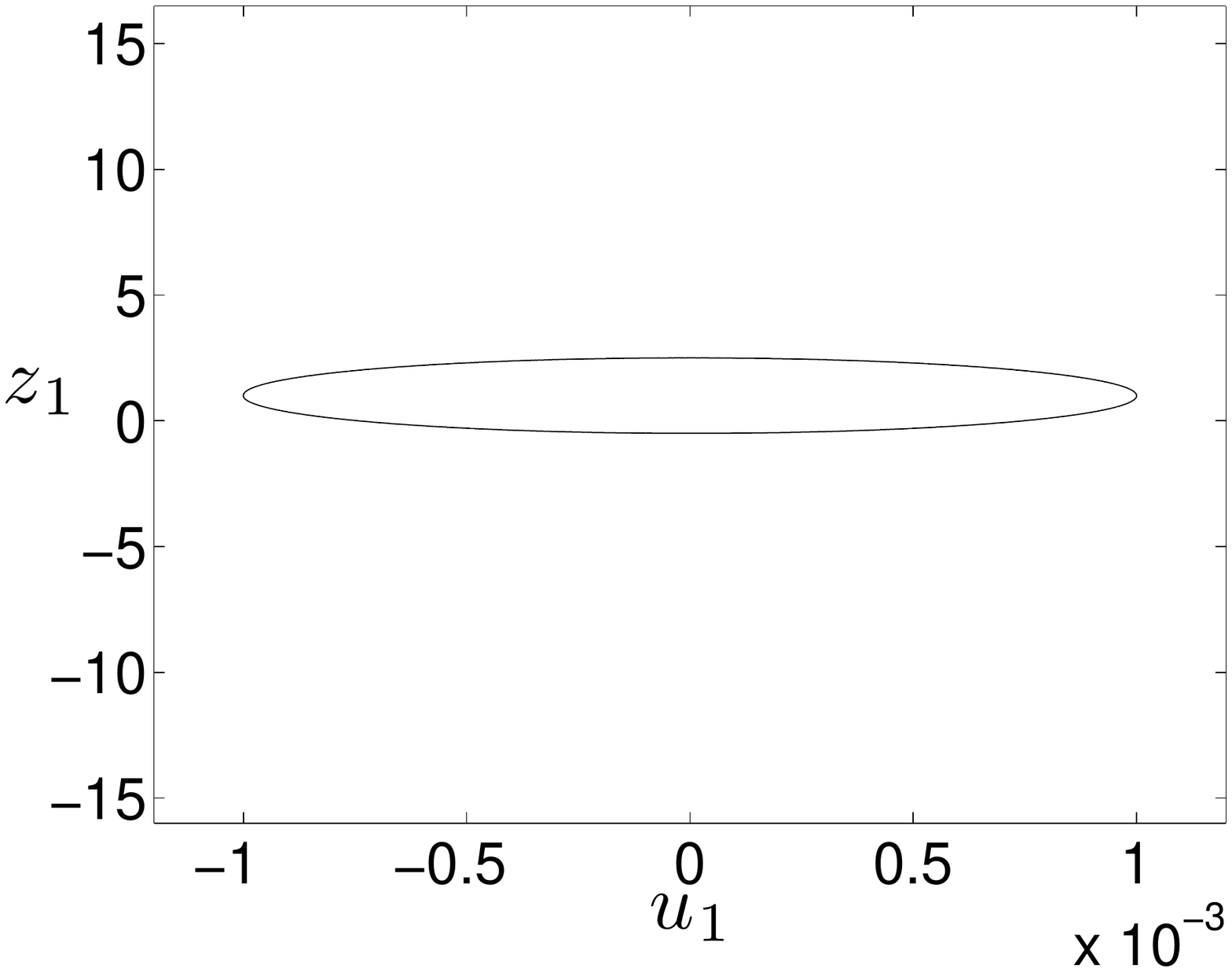}} 
    \subfloat[$\omega=0.001$]{ \includegraphics[trim = 0mm 60mm 0mm 60mm, clip,width=0.26\textwidth]{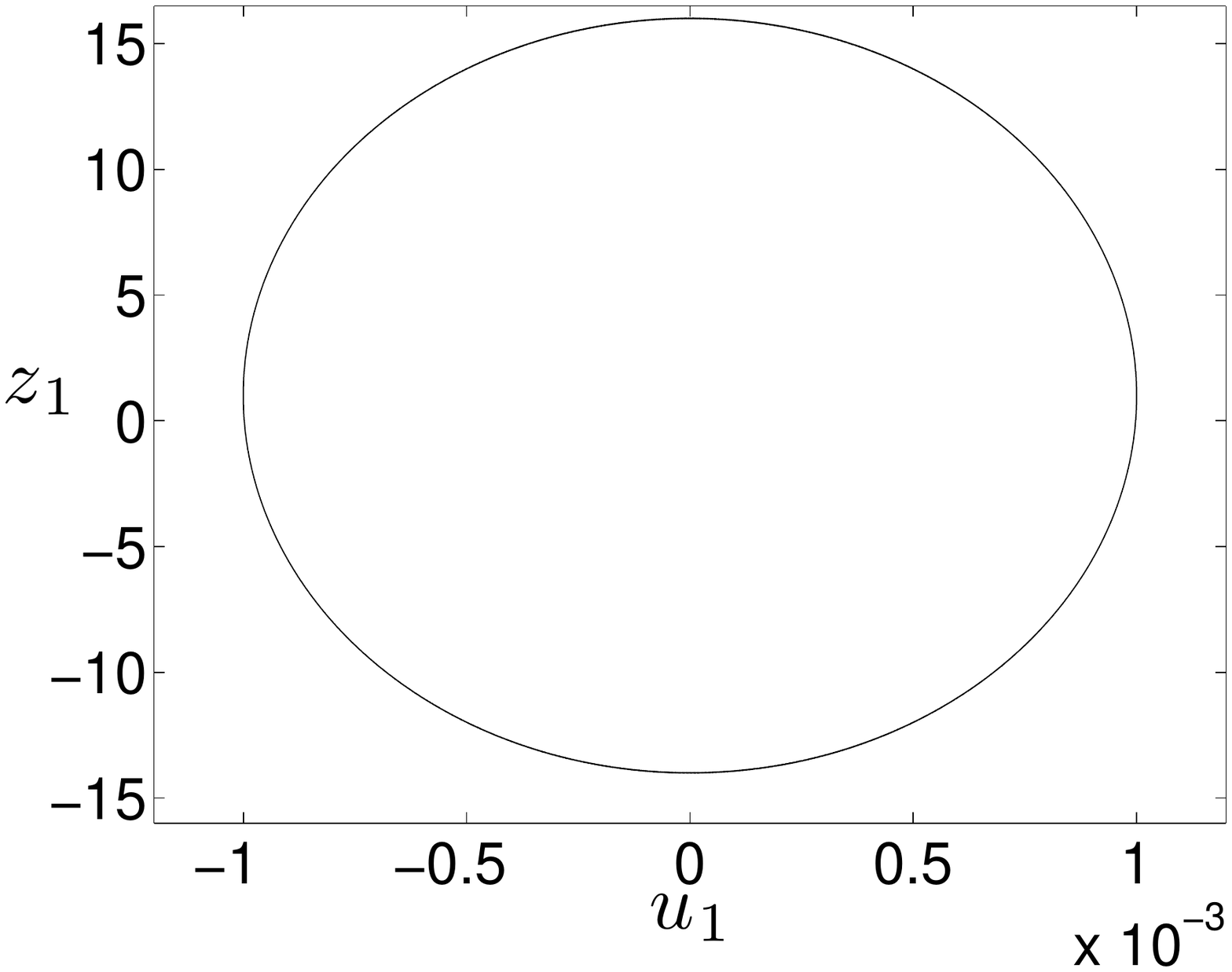}}} 
\caption{\label{figlinearhysteresisloopM1} \small Hysteresis loops for $z_1(x,t)$ of the linear Landau-Lifshitz equation with $x=0.6$ and $\nu=0.02$.  The linearization is at $\mathbf a=(1,0,0)$. The input is $\mathbf u(t)=\left(0.001\cos(\omega t),0,0\right)$ and the initial condition is $\mathbf z_0(x)=\left( 1, 0,0 \right)$.}
\end{figure}

Using Theorem~\ref{thmzeroeigenvalue} and Definition~\ref{defmultiequilibrium} indicates the linear Landau-Lifshitz equation exhibits hysteresis. Furthermore, simulations with   periodic inputs  were performed to determine whether persistent loops exist and hence Definition~\ref{defBernstein} is satisfied.  Again, $\nu=0.02$,  $L=1$,  and the same periodic input is applied as for the nonlinear Landau-Lifshitz equation.  Figure~\ref{figlinearhysteresisloopM1} shows the input-output curves for the first component of the solution to the linear Landau-Lifshitz equation with $\mathbf a=(1,0,0)$ and initial condition  $\mathbf z_0(x)=\left( 1, 0,0 \right)$.   From the figure, it is clear a loop persists  as the frequency of the input approaches zero. Similar plots are obtained when the control is on the second and third components. The hysteresis loops in Figure~\ref{figlinearhysteresisloopM1}  are similar in shape to  the nonlinear Landau-Lifshitz equation depicted in Figure~\ref{figLoopFrequencyVaried}, both of which have a continuum of equilibria.

\section{Controller Design}\label{seccontrol}

A control, $ \mathbf u(t)$, is introduced into the Landau-Lifshitz equation~(\ref{eqLLGuoDing}) as follows
\begin{align}
\label{eqcontrolledLL}
\frac{\partial \mathbf m}{\partial t} & = \mathbf m \times  \mathbf m_{xx}-\nu\mathbf m\times\left(\mathbf m\times \mathbf m_{xx}\right)+ \mathbf u(t)\\
 \mathbf m(x,0)&=\mathbf m_0(x)\nonumber\\
 \mathbf m_x(0,t)&=\mathbf m_x(L,t)=\mathbf 0. \nonumber
\end{align}
The goal is to choose a control $\mathbf u(t)$  so the system governed by the Landau-Lifshitz equation moves from an arbitrary initial condition, possibly an equilibrium point, to a specified equilibrium point $\mathbf r$, where $\mathbf r \in E$ and $E$ is defined in (\ref{equilibriumset}).  The control function needs to be chosen so that  $\mathbf r $  becomes an asymptotically stable equilibrium point of the controlled system.  

\begin{figure}[h]
\scalebox{1}{\begin{picture}(450,140)(-120,-120)
\put(-30,0){\vector(1,0){50}}
\put(-35,0){\circle{10}}
\put(-70,0){\vector(1,0){30}}
\put(20,-20){\framebox(20,40)[c]{ $k$ }}
\put(40,0){\vector(1,0){70}}
\put(110,-20){\framebox(140,40)[c]{Landau-Lifshitz Equation}}
\put(250,0){\vector(1,0){40}}
\put(270,0){\line(0,-1){40}}
\put(270,-40){\line(-1,0){305}}
\put(-35,-40){\vector(0,1){35}}
\put(-45,-10){-}
\put(-90,-2){\small{ $\mathbf r$}}
\put(60,6){\small{ $\mathbf u(x,t)$}}
\put(295,-2){\small{$\mathbf m(x,t)$}}
\end{picture}}\vspace*{-2.25cm}
\caption{\label{figblockdiagram} \small  Closed-loop control of  the Landau-Lifshitz equation.}\end{figure}
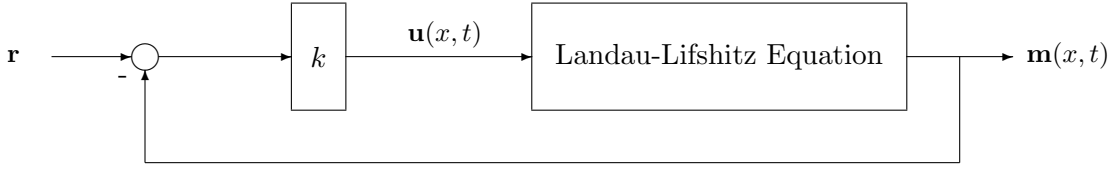

 Theorem~\ref{thmzeroeigenvalue} implies zero is an eigenvalue of the   uncontrolled linearized Landau-Lifshitz equation. For finite-dimensional linear systems,  simple proportional control  of a system with a zero eigenvalue yields asymptotic tracking of a specified  state and this motivates choosing the control
\begin{equation}\label{equ}
\mathbf u(x,t)=k(\mathbf r -\mathbf m(x,t) )
\end{equation}
where $k$ is a positive constant control parameter for equation~(\ref{eqcontrolledLL}). It is clear that $\mathbf r$ is an equilibrium point of (\ref{eqcontrolledLL}) with the control in (\ref{equ}). Figure~\ref{figblockdiagram} is a block diagram representation of (\ref{eqcontrolledLL}) with control  (\ref{equ}).

The following theorem establishes well-posedness of the controlled equation. In particular, for any initial condition $\mathbf m_0,$ the solution to (\ref{eqcontrolledLL}) with control  $u(t) =k(\mathbf r-{\mathbf m}) $ satisfies
$\| \mathbf m(\cdot, t)\|_{\mathcal L_2^3} \leq 1.$

\begin{thm}\label{thmstrongsolutionLLcontrol}
For any $\mathbf r \in E,$ define the operator
\begin{equation}
B\mathbf m =k(\mathbf r-{\mathbf m}) \label{eqlinearAforcontrolledLL}.
\end{equation}
If $k>0$, the nonlinear operator $f+B$ with domain $D$, where $f$ and  $D$ are defined in (\ref{defn:f}), (\ref{setDforfullLL}) respectively, generates a nonlinear contraction semigroup on  ${\mathcal L_2^3}$.
\end{thm}

\begin{pf}

(i) For any $\mathbf m, \mathbf y \in D$,
\begin{align*}
&\langle f(\mathbf m)+B\mathbf m -  (f(\mathbf y) +B\mathbf y),\mathbf m -\mathbf y \rangle_{\mathcal L_2^3}  \\
&=\langle f(\mathbf m)-f(\mathbf y), \mathbf m -\mathbf y \rangle_{\mathcal L_2^3}  +\langle B\mathbf m-B\mathbf y,\mathbf m -\mathbf y \rangle_{\mathcal L_2^3}.
\end{align*}
Since $f$ generates a nonlinear contraction semigroup (Theorem~\ref{thmuncontrolledLLsemigroup}), then $f$ is dissipative \cite[Proposition~2.98]{Luo1999}; that is,
\[
\langle f(\mathbf m)-f(\mathbf y), \mathbf m -\mathbf y \rangle_{\mathcal L_2^3}\leq 0.
\]
It follows that 
\begin{align*}
\langle f(\mathbf m)+B\mathbf m -  (f(\mathbf y) +B\mathbf y),\mathbf m -\mathbf y \rangle_{\mathcal L_2^3}  
&\leq \langle B\mathbf m-B\mathbf y,\mathbf m -\mathbf y \rangle_{\mathcal L_2^3}\\
&=\langle -k\mathbf m+k\mathbf y,\mathbf m -\mathbf y \rangle_{\mathcal L_2^3}\\
&=-k\langle \mathbf m-\mathbf y,\mathbf m -\mathbf y \rangle_{\mathcal L_2^3}\\
&\leq 0
\end{align*}
and hence $f+B$ is dissipative.

(ii) Since $f$ generates a nonlinear contraction semigroup (Theorem~\ref{thmuncontrolledLLsemigroup}),  $\mathrm{ran}(\mathrm{I}-\hat \alpha f)=\mathcal L_2^3$ for any $\hat\alpha >0$  \cite[Lemma~2.1]{Kato1967}.  
This means  that for any $\mathbf y_2 \in \mathcal L_2^3$ there exists $\mathbf m \in D$ such that $\mathbf{m}-\hat{\alpha}f(\mathbf m)=\mathbf y_2. $ 
Choose any $\mathbf y_1 \in \mathcal L_2^3$, $\alpha >0$ and define 
\[
\mathbf y_2 = \frac{\mathbf y_1}{1+\alpha k}+\frac{\alpha k\mathbf r }{1+\alpha k}
\]
and
\[
\hat \alpha=\frac{\alpha  }{1+\alpha k}.
\]
There exists $\mathbf m \in D$ such that
\[
\mathbf{m}-\frac{\alpha }{1+\alpha k}f(\mathbf m)=\mathbf y_2=  \frac{\mathbf y_1}{1+\alpha k}+\frac{\alpha k\mathbf r }{1+\alpha k}.
\]
Solving for $\mathbf y_1$ leads to
\[
\mathbf y_1 = \mathbf m-\alpha (k(\mathbf r-\mathbf m)+ f(\mathbf m)).
\]
Thus, for any $\mathbf y_1 \in \mathcal L_2^3$, there exists $\mathbf m \in D$ such that $
\mathbf y_1 = (\mathrm I-\alpha (B+f))\mathbf m$ and hence $\mathrm{ran}\left(\mathrm I-\alpha (B+f)\right)=\mathcal L_2^3$ for some $\alpha >0$. It follows that the range is $\mathcal L_2^3$  for all  $\alpha>0$ \cite[Lemma~2.1]{Kato1967}. 

Thus, since $B+f$ is dissipative and the range of $(\mathrm I-\alpha (B+f))$ is $\mathcal L_2^3$, then $B+f$ generates a nonlinear contraction semigroup \cite[Proposition~2.114]{Luo1999}.
 \qed
 \end{pf}

The following lemmas are needed in the proof of the main results in Theorem~\ref{thmr0isasymstable} and Theorem~\ref{thmr0isexpstable}.  The first theorem demonstrates the control in (\ref{equ}) can steer the dynamics in the Landau-Lifshitz to an asymptotically stable state in the $L_2^3$-norm, while the latter theorem establishes exponential stability in the $H_1$-norm, $||\mathbf m ||^2_{H_1}=||\mathbf m||_{\mathcal L_2^3}^2+||\mathbf m_x||_{\mathcal L_2^3}^2.$  

\begin{lem}\label{lemmaboundonacrossz}
If $\mathbf a \in E$ where $E$ is defined in (\ref{equilibriumset}), then $||\mathbf a \times \mathbf  m_{}||_{\mathcal L_2^3} \leq || \mathbf m_{}||_{\mathcal L_2^3}$ for all $\mathbf m \in \mathcal L_2^3$. 
\end{lem}
\begin{pf}
Since $||\mathbf a \times \mathbf  m||_2 =||\mathbf a||_2||\mathbf m||_2\sin(\theta)$ where $\theta$ is the angle between $\mathbf a$ and $\mathbf m$, and $||\mathbf a||_2=1$, then
$||\mathbf a \times \mathbf  m||_2 \leq ||\mathbf m||_2$. Extending to the $L_2^3$-norm, the desired result is obtained.
\end{pf}

Lemmas ~\ref{thmderivativemcrossmprime} and~\ref{thmderivativemcrossmprimedotproduct} are simple consequences of the product rule.
\begin{lem}\label{thmderivativemcrossmprime}
For $\mathbf m\in \mathcal L_2^3$, the derivative of $\mathbf g=\mathbf m \times \mathbf m_x$ is $\mathbf g_x=\mathbf m \times \mathbf m_{xx}$.
\end{lem}

\begin{lem}\label{thmderivativemcrossmprimedotproduct}
For $\mathbf m \in \mathcal L_2^3$, the derivative of $f=\left(\mathbf m \times \mathbf m_x\right)^{\mathrm T}\left(\mathbf m \times \mathbf m_x\right)$ is $f_x=2\left(\mathbf m \times \mathbf m_{x}\right)^{\mathrm T}\left(\mathbf m \times \mathbf m_{xx}\right).$
\end{lem}

\begin{lem}
\label{lemmazerointegral}
For $\mathbf m \in \mathcal L_2^3$ satisfying (\ref{eqboundarycondition}), 
\[
\int_0^L (\mathbf m-\mathbf r)^{\mathrm T}(\mathbf m \times \mathbf m_{xx})dx=0.
\]
\end{lem}
\begin{pf}
Integrating by parts, and applying  Lemma~\ref{thmderivativemcrossmprime} and the boundary conditions (\ref{eqboundarycondition}) implies 
\[
\int_0^L (\mathbf m-\mathbf r)^{\mathrm T}(\mathbf m \times \mathbf m_{xx})dx = -\int_0^L \mathbf m_x^{\mathrm T}(\mathbf m \times \mathbf m_{x})dx.
\]
From properties of cross products,
$
\mathbf m_x^{\mathrm T}(\mathbf m \times \mathbf m_{x}) =  \mathbf m^{\mathrm T}(\mathbf m_x \times \mathbf m_{x})=0,
$
and hence the integral is zero. \qed
\end{pf}
\begin{lem}\label{lemmaPoincareInequalityforCrossProducts}
For $\mathbf m \in \mathcal L_2^3$ satisfying (\ref{eqboundarycondition}), 
\[
|| \mathbf m \times \mathbf m_x ||_{\mathcal L_2^3} \leq 4L^2|| \mathbf m \times \mathbf m_{xx} ||_{\mathcal L_2^3}
\]
\end{lem}
\begin{pf}
Integrating by parts, using Lemma~\ref{thmderivativemcrossmprimedotproduct} and the boundary conditions (\ref{eqboundarycondition}) leads to
\[
|| \mathbf m \times \mathbf m_x ||_{\mathcal L_2^3}^2 =-\int_0^L2\left(\mathbf m \times \mathbf m_{x}\right)^{\mathrm T}\left(\mathbf m \times \mathbf m_{xx}\right)xdx.
\]
 It follows from Young's inequality that
\begin{align*}
|| \mathbf m \times \mathbf m_x ||_{\mathcal L_2^3}^2 
 \leq \frac{1}{2}\int_0^L\left(\mathbf m \times \mathbf m_{x}\right)^{\mathrm T}\left(\mathbf m \times \mathbf m_{x}\right)dx+\int_0^L2\left(\mathbf m \times \mathbf m_{xx}\right)^{\mathrm T}\left(\mathbf m \times \mathbf m_{xx}\right)x^2dx .
 \end{align*}
 Since $x \in [0,L],$ 
 \begin{align*} 
 || \mathbf m \times \mathbf m_x ||_{\mathcal L_2^3}^2
 \leq \frac{1}{2}||\mathbf m \times \mathbf m_{x}||_{\mathcal L_2^3}^{2}+2L^2||\mathbf m \times \mathbf m_{xx}||_{\mathcal L_2^3}^2.
\end{align*}
Rearranging gives the desired inequality. \qed
\end{pf}

\begin{thm}\label{thmr0isasymstable}
Let $\mathbf r$ be an equilibrium point of (\ref{eqcontrolledLL}) with control defined in (\ref{equ}). For any positive constant $k$ such that $k> 8\nu L^4$, $\mathbf r$ is a globally asymptotically stable point of (\ref{eqcontrolledLL})  in the $\mathcal L_2^3$--norm.
\end{thm}

\begin{pf}
The Lyapunov candidate is
\[
V(\mathbf m)=\frac{1}{2}\left| \left| \mathbf m-\mathbf r\right|\right|_{\mathcal L_2^3}^2+\frac{1}{2}\left| \left|  \mathbf m_x\right|\right|_{\mathcal L_2^3}^2
\]
which is clearly nonegative. Furthermore, $V=0$ if and only if $\mathbf m=\mathbf r$.
Taking the derivative of $V$
\begin{align*}
\frac{dV}{dt}&=\int_0^L(\mathbf m -\mathbf r)^{\mathrm T}\dot{{\mathbf m}} dx+\int_0^L\mathbf m_x^{\mathrm T} \dot{{\mathbf m}}_xdx \\
&=\int_0^L(\mathbf m -\mathbf r)^{\mathrm T}\dot{\mathbf m} dx-\int_0^L\mathbf m_{xx}^{\mathrm T} \dot{\mathbf m}dx
\end{align*}
where the dot notation means differentiation with respect to $t$. Substituting in  (\ref{eqcontrolledLL}) to eliminate $\dot{\mathbf m}$,
\begin{align*}
\frac{dV}{dt}
&= \int_0^L(\mathbf m -\mathbf r)^{\mathrm T}  \left(\mathbf m \times  \mathbf m_{xx}\right)dx-\nu\int_0^L (\mathbf m -\mathbf r)^{\mathrm T}\left(\mathbf m\times\left(\mathbf m\times \mathbf m_{xx}\right)\right)dx \\
&+k\int_0^L(\mathbf m -\mathbf r)^{\mathrm T} (\mathbf r -\mathbf m)dx -\int_0^L\mathbf m_{xx}^{\mathrm T}  \left(\mathbf m \times  \mathbf m_{xx}\right)dx\\
&+\nu\int_0^L \mathbf m_{xx}^{\mathrm T}\left(\mathbf m\times\left(\mathbf m\times \mathbf m_{xx}\right)\right)dx -k\int_0^L\mathbf m_{xx}^{\mathrm T} (\mathbf r -\mathbf m)dx.
\end{align*}
From Lemma~\ref{lemmazerointegral}, the first integral is zero. Furthermore, from properties of cross products,
\[
\mathbf m_{xx}^{\mathrm T}  \left(\mathbf m \times  \mathbf m_{xx}\right)=\mathbf m^{\mathrm T}  \left(\mathbf m_{xx} \times  \mathbf m_{xx}\right)=0,
\] 
and hence
\[
\int_0^L\mathbf m_{xx}^{\mathrm T}  \left(\mathbf m \times  \mathbf m_{xx}\right)dx=0.
\]
It follows that 
\begin{align*}
\frac{dV}{dt}=&-\nu\int_0^L (\mathbf m -\mathbf r)^{\mathrm T}\left(\mathbf m\times\left(\mathbf m\times \mathbf m_{xx}\right)\right)dx -k||\mathbf m -\mathbf r||_{\mathcal L_2^3}^2\\
&-\nu||\mathbf m\times \mathbf m_{xx}||_{\mathcal L_2^3}^2 -k\int_0^L\mathbf m_{xx}^{\mathrm T} (\mathbf r -\mathbf m)dx.
\end{align*}
Applying integration by parts to the last integral leads to
\begin{align}
\frac{dV}{dt}&=-\nu\int_0^L (\mathbf m -\mathbf r)^{\mathrm T}\left(\mathbf m\times\left(\mathbf m\times \mathbf m_{xx}\right)\right)dx -k||\mathbf m -\mathbf r||_{\mathcal L_2^3}^2-\nu||\mathbf m\times \mathbf m_{xx}||_{\mathcal L_2^3}^2 -k||\mathbf m_{x}||_{\mathcal L_2^3}^2 \nonumber\\
&=-\nu\int_0^L\left( (\mathbf m -\mathbf r)\times \mathbf m\right)^{\mathrm T}\left(\mathbf m\times \mathbf m_{xx}\right)dx -k||\mathbf m -\mathbf r||_{\mathcal L_2^3}^2-\nu||\mathbf m\times \mathbf m_{xx}||_{\mathcal L_2^3}^2 -k||\mathbf m_{x}||_{\mathcal L_2^3}^2 \nonumber\\
&=\nu\int_0^L\left( \mathbf r\times \mathbf m\right)^{\mathrm T}\left(\mathbf m\times \mathbf m_{xx}\right)dx -k||\mathbf m -\mathbf r||_{\mathcal L_2^3}^2-\nu||\mathbf m\times \mathbf m_{xx}||_{\mathcal L_2^3}^2 -k||\mathbf m_{x}||_{\mathcal L_2^3}^2. \label{eqderivativeofVlastintegral}
\end{align}
Applying integration by parts with Lemma~\ref{thmderivativemcrossmprime} to the integral implies
\begin{align*}
\int_0^L\left( \mathbf r\times \mathbf m\right)^{\mathrm T}\left(\mathbf m\times \mathbf m_{xx}\right)dx&=\left[\left(\mathbf r \times \mathbf m\right)^{\mathrm T}\left(\mathbf m \times \mathbf m_x\right)\right]_0^L-\int_0^L\left( \mathbf r\times \mathbf m_x\right)^{\mathrm T}\left(\mathbf m\times \mathbf m_{x}\right)dx
\end{align*}
and substituting in the boundary conditions in (\ref{eqboundarycondition}) leads to 
\begin{align*}
\int_0^L\left( \mathbf r\times \mathbf m\right)^{\mathrm T}\left(\mathbf m\times \mathbf m_{xx}\right)dx&=-\langle  \mathbf r\times \mathbf m_x, \mathbf m\times \mathbf m_{x}\rangle_{\mathcal L_2^3}.
\end{align*}
Then from Cauchy-Schwarz and Lemma~\ref{lemmaPoincareInequalityforCrossProducts},
\begin{align*}
\int_0^L\left( \mathbf r\times \mathbf m\right)^{\mathrm T}\left(\mathbf m\times \mathbf m_{xx}\right)dx&\leq|| \mathbf r\times \mathbf m_x||_{\mathcal L_2^3} ||\mathbf m\times \mathbf m_{x}||_{\mathcal L_2^3}
\\
&\leq 4L^2  || \mathbf r\times \mathbf m_x||_{\mathcal L_2^3} ||\mathbf m\times \mathbf m_{xx}||_{\mathcal L_2^3}.
\end{align*}
It follows from Young's Inequality that
\begin{align*}
\int_0^L\left( \mathbf r\times \mathbf m\right)^{\mathrm T}\left(\mathbf m\times \mathbf m_{xx}\right)dx&\leq 8L^4 || \mathbf r\times \mathbf m_x||_{\mathcal L_2^3} ^2+\frac{1}{2}||\mathbf m\times \mathbf m_{xx}||_{\mathcal L_2^3}^2
\end{align*}
and from Lemma~\ref{lemmaboundonacrossz},
\begin{align*}
\int_0^L\left( \mathbf r\times \mathbf m\right)^{\mathrm T}\left(\mathbf m\times \mathbf m_{xx}\right)dx&\leq 8L^4 ||  \mathbf m_x||_{\mathcal L_2^3} ^2+\frac{1}{2} ||\mathbf m\times \mathbf m_{xx}||_{\mathcal L_2^3}^2
\end{align*}
Substituting this result into (\ref{eqderivativeofVlastintegral}) leads to
\begin{equation}
\label{eqderivativeofVLL}
\frac{dV}{dt}\leq -\left(k-8\nu L^4\right)||\mathbf m_{x}||_{\mathcal L_2^3}^2-\frac{\nu}{2}||\mathbf m\times \mathbf m_{xx}||_{\mathcal L_2^3}^2-k||\mathbf m -\mathbf r||_{\mathcal L_2^3}^2.
\end{equation}
The derivative is negative if $k > 8 \nu L^4.$ 
It follows that 
\[
\frac{dV}{dt}\leq-k||\mathbf m -\mathbf r||_{\mathcal L_2^3}^2.
\]
Therefore, $dV/dt <0$ for all $\mathbf m \neq \mathbf r$ and $dV/dt=0 $ if $\mathbf m =\mathbf r$. Since $V(\mathbf m)\geq \frac{1}{2}||\mathbf m -\mathbf r||_{\mathcal L_2^3}^2,$ $V\rightarrow \infty$ as $||\mathbf m -\mathbf r||\rightarrow \infty$.
 From Lyapunov's Theorem \cite[Theorem~6.2.13]{Michel1995}, $\mathbf r$ is a globally asymptotically stable equilibrium of (\ref{eqcontrolledLL}).  \qed
 \end{pf}

\begin{thm}\label{thmr0isexpstable} 
Let $\mathbf r$ be an equilibrium point of (\ref{eqcontrolledLL}) with control defined in (\ref{equ}).  For any positive constant $k$ such that ${k> 8\nu L^4}$, $\mathbf r$ is  a globally exponentially stable equilibrium point of (\ref{eqcontrolledLL}) in the $H_1$--norm. That is, for any initial condition on $H_1$, $\mathbf m$ decreases exponentially in the $H_1$-norm to $\mathbf r.$
\end{thm}

\begin{pf} 
In the proof of Theorem~\ref{thmr0isasymstable}, we have from (\ref{eqderivativeofVLL}) that the derivative of V satisfies
\begin{align*}
\frac{dV}{dt}&\leq -\left(k- {8\nu L^4}\right)||\mathbf m_{x}||_{\mathcal L_2^3}^2-\frac{\nu}{2}||\mathbf m\times \mathbf m_{xx}||_{\mathcal L_2^3}^2-k||\mathbf m -\mathbf r||_{\mathcal L_2^3}^2
\end{align*}
and hence
\begin{align*}
\frac{dV}{dt}&\leq -\left(k-{8\nu L^4}\right)||\mathbf m_{x}||_{\mathcal L_2^3}^2-k||\mathbf m -\mathbf r||_{\mathcal L_2^3}^2\\
&\leq -\left(k-{8\nu L^4}\right)\left(||\mathbf m_{x}||_{\mathcal L_2^3}^2+||\mathbf m -\mathbf r||_{\mathcal L_2^3}^2\right)\\
&=-2\left(k-{8\nu L^4}\right)V.
\end{align*}
Integrating with respect to time 
\begin{align*}
||\mathbf m_{x}||_{\mathcal L_2^3}^2+||\mathbf m -\mathbf r||_{\mathcal L_2^3}^2
& \leq e^{-2\left(k-{8\nu L^4}\right)t}\left(||\mathbf {m}_{x}(x,0)||_{\mathcal L_2^3}^2+||\mathbf m(x,0) -\mathbf r||_{\mathcal L_2^3}^2 \right)
\end{align*}
Noting that $\mathbf r$ does not depend on $x$, it follows that
\[
||\mathbf m-\mathbf r||_{H_1}^2\leq e^{-2\left(k-{8\nu L^4}\right)t}||\mathbf m(x,0) -\mathbf r||_{H_1}^2
\]
and since $k>{8\nu L^4}$, $\mathbf r$ is an exponentially stable equilibrium point of (\ref{eqcontrolledLL}). \qed
\end{pf}

A natural question is whether  $\mathbf r$ is exponentially stable in the $\mathcal L_2^3$-norm. 
Analysis of the linear Landau-Lifshitz equation provides insight to this question. For the control in (\ref{equ}), the linearized controlled Landau-Lifshitz equation is 
\begin{align}\label{eqcontrolledlinearLL}
\frac{\partial \mathbf z}{\partial t} & =\nu \mathbf  z_{xx}+\mathbf a \times  {\mathbf z}_{xx} +k \left(\mathbf r-\mathbf z\right), \qquad \mathbf z(0)=\mathbf z_0
\end{align}
with the same  boundary conditions $\mathbf z_{x}(0)=\mathbf z_{x}(L)=\mathbf 0 .$ 
 Since the uncontrolled linear Landau-Lifshitz equation (\ref{eqstatespaceform})  generates a linear semigroup and $k \left(\mathbf r-\mathbf z\right)$ is a bounded linear (affine) operator, then the operator in (\ref{eqcontrolledlinearLL}) generates a semigroup \cite[Theorem~3.2.1]{Curtain1995}.  Substituting $\mathbf z=\mathbf{r}$ into (\ref{eqcontrolledlinearLL}) leads to ${\partial \mathbf z}/{\partial t}=\mathbf 0$ and hence $\mathbf{r}$ is a stable  equilibrium point of (\ref{eqcontrolledlinearLL}). 
  
\begin{thm}\label{thmlinearr0isasymstable}
Let $\mathbf r \in E .$  For any positive constant $k$, $\mathbf r$ is an exponentially stable equilibrium of the linearized system (\ref{eqcontrolledlinearLL})  in
the $\mathcal L_2^3$--norm. 
\end{thm}

\begin{pf}
For $\mathbf z \in D(A)$, where $D(A)=D$ as in equation~(\ref{setDforfullLL}), consider the Lyapunov candidate
\begin{equation*}
V(\mathbf z)=\frac{1}{2}\left| \left|\mathbf z-\mathbf r\right|\right|_{\mathcal L_2^3}^2.
\end{equation*}
It is clear that $V\geq 0$ for all $\mathbf z\in D(A)$ and furthermore, $V(\mathbf z)=0$ only when $\mathbf z=\mathbf r$.  Therefore, $V(\mathbf z)>0$ for all $\mathbf z \in D(A)\backslash  \{\mathbf r\}$.

Taking the derivative of $V(\mathbf z)$ implies
\begin{align*}
\frac{dV}{dt}&=\int_0^L(\mathbf z-\mathbf r)^{\mathrm T}\dot{ {\mathbf z}} dx.
\end{align*}
Substituting in (\ref{eqcontrolledlinearLL}) yields 
\begin{align*}
\frac{dV}{dt}&=\nu\int_0^L(\mathbf z-\mathbf r)^{\mathrm T} {\mathbf z}_{xx}dx +\int_0^L(\mathbf r-\mathbf z)^{\mathrm T}\left({\mathbf a \times {\mathbf  z_{xx}}}\right)dx +k\int_0^L(\mathbf z-\mathbf r)^{\mathrm T} ({\mathbf r-\mathbf z}) dx.
\end{align*}
By Lemma~\ref{lemmazerointegral}, the middle term is zero. Using integration by parts,  the first term becomes
\[
-\nu\int_0^L\mathbf z_x^{\mathrm T} {\mathbf z}_{x}dx.
\]
It follows that 
\begin{align*}
\frac{dV}{dt}&=-\nu||\mathbf z_x||_{\mathcal L_2^3}^2 -k||\mathbf z-\mathbf r||_{\mathcal L_2^3}^2
\end{align*}
and since $\nu\geq0$, 
\begin{align*}
\frac{dV}{dt}&\leq -k||\mathbf z-\mathbf r||_{\mathcal L_2^3}^2=-2kV.
\end{align*}
Solving yields
\[
||\mathbf z-\mathbf r||_{\mathcal L_2^3}^2 \leq  e^{-2kt}||\mathbf z_0-\mathbf r||_{\mathcal L_2^3}^2.
\]
For $k>0$ the equilibrium point, $\mathbf r,$ of (\ref{eqcontrolledlinearLL}) is exponentially stable. \qed
\end{pf}
Theorem~\ref{thmlinearr0isasymstable} suggests the equilibrium point in the controlled nonlinear Landau-Lifshitz equation (\ref{eqcontrolledLL}) is exponentially stable in the $\mathcal L_2^3$--norm. However, since the nonlinearity in the Landau-Lifshitz equation is unbounded, stability of the linear equation does not necessarily reflect stability of the original nonlinear equation; see  \cite{CoronNguyen2015,alJamal2013,AM2014}.

\section{Example}\label{secExample}
\begin{figure}[h]
\centering\scalebox{0.7}{
\includegraphics{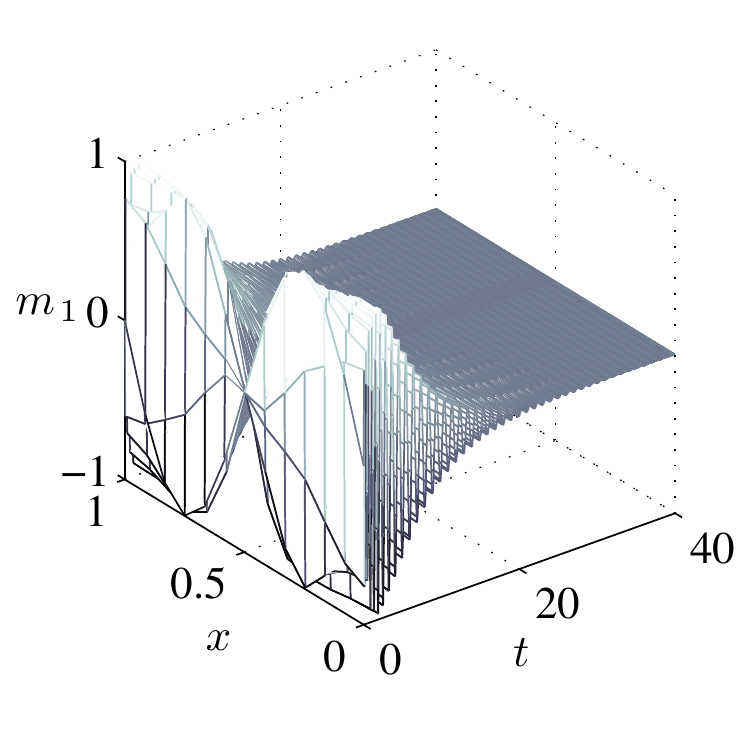}
\includegraphics{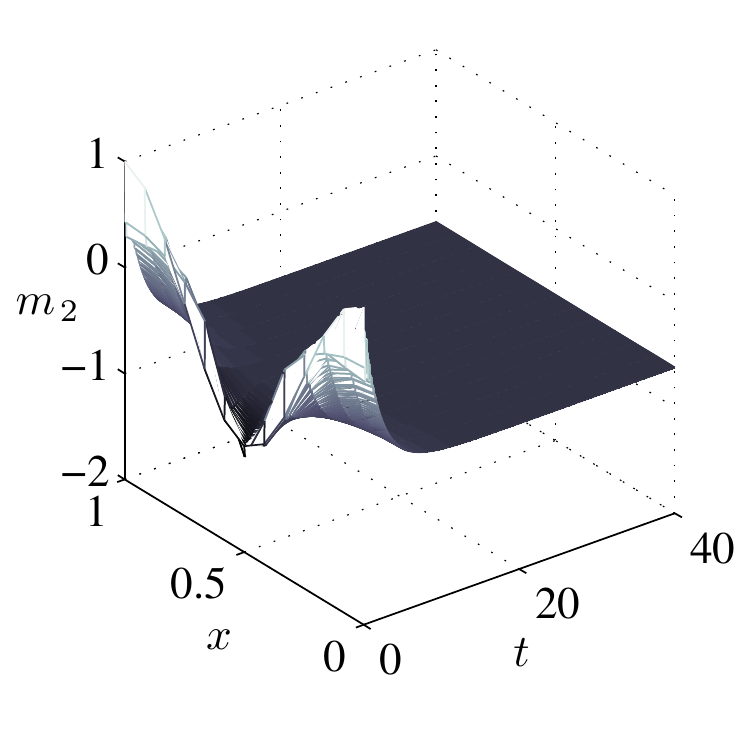}
\includegraphics{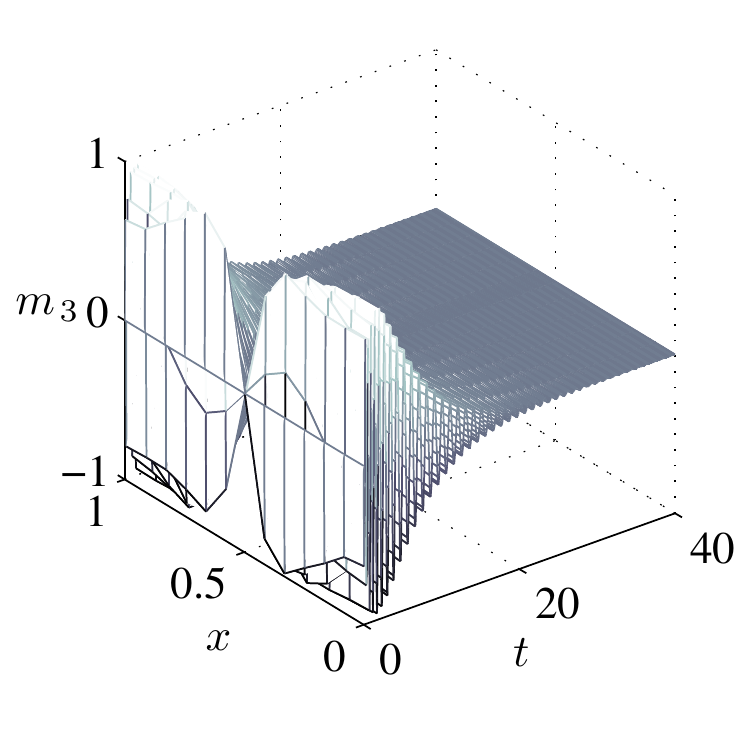}}
 \caption{\label{figMagnetizationNoControl3d}
   Magnetization in the uncontrolled  (nonlinear) Landau-Lifshitz equation  moves from initial condition  $\mathbf{{m}_0} (x)$,  to the  equilibrium $\mathbf {r_0}=(0,-0.6,0)$.   }
\end{figure}

Simulations illustrating the stabilization of the (nonlinear) Landau-Lifshitz equation are constructed using a Galerkin approximation with 12 linear spline elements. For the following simulations, the parameters  are $\nu=0.02$ and $L=1$ with  initial condition  $\mathbf {{m}_0} (x)=(\sin(2\pi x),\cos(2\pi x),0)$. Figure~\ref{figMagnetizationNoControl3d} illustrates the solution to the uncontrolled Landau-Lifshitz equation settles to $\mathbf{r_0}=(0,-0.6,0).$ 

\begin{figure}[h]
\centering\scalebox{0.7}{
\includegraphics{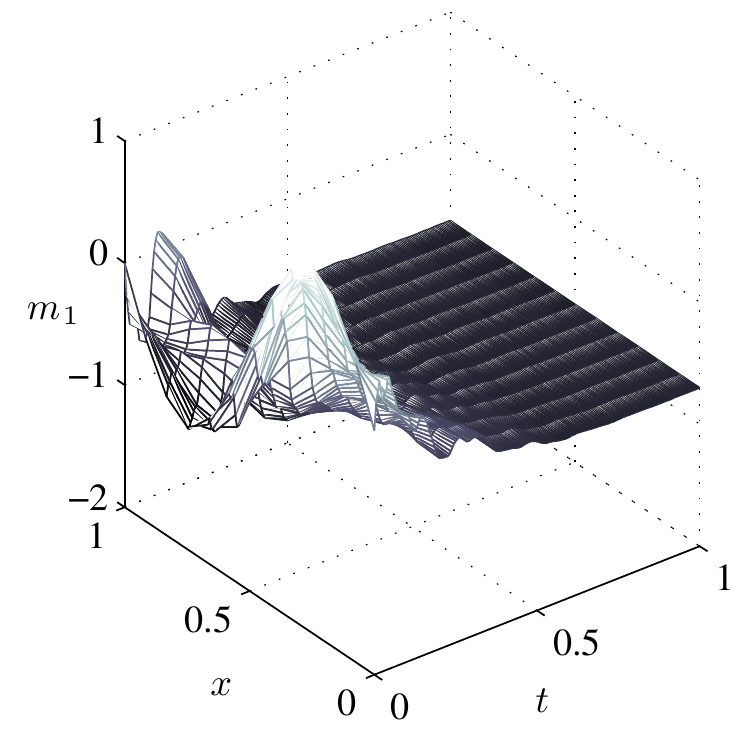}
\includegraphics{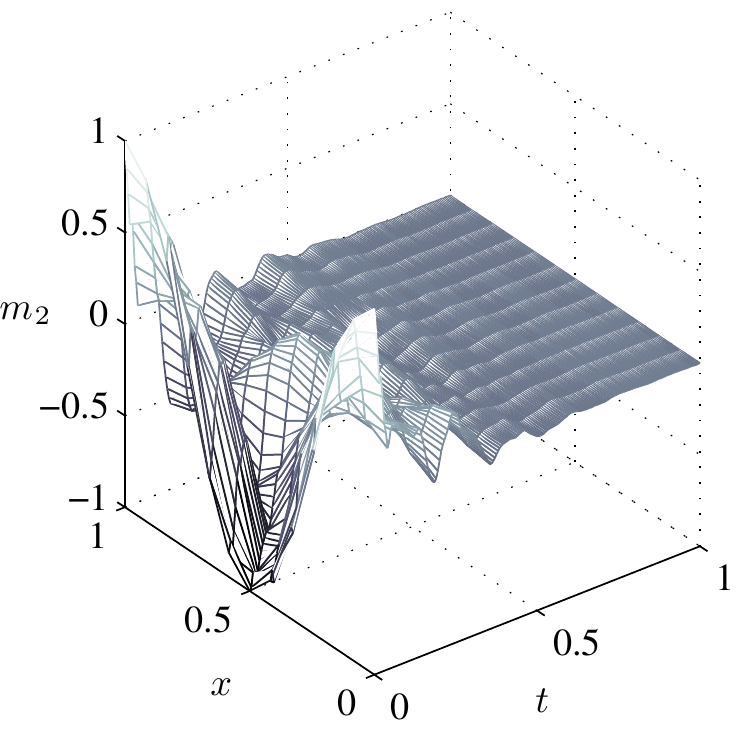}
\includegraphics{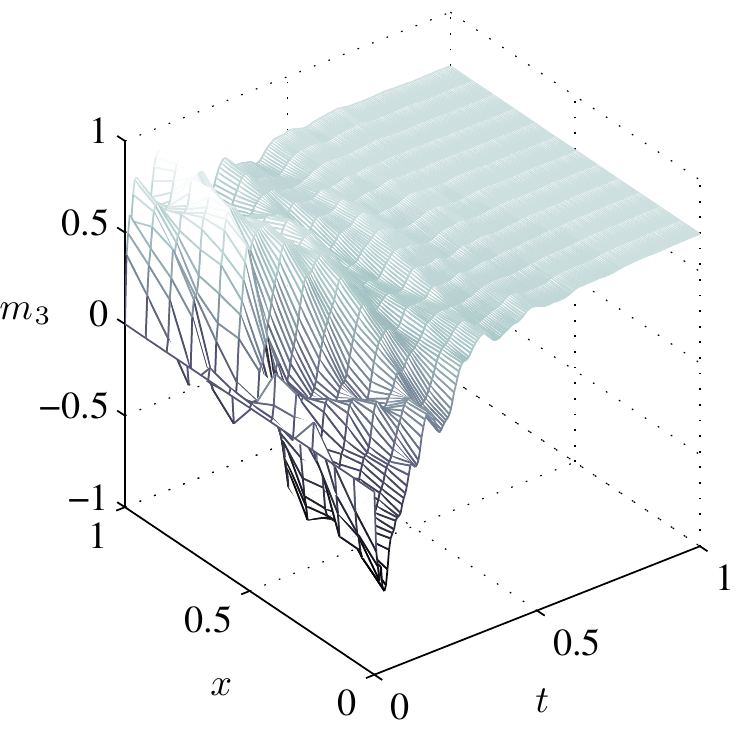}}
 \caption{\label{figMagnetization3d} With a proportional control $(k=0.5)$,  magnetization in the (nonlinear) Landau-Lifshitz equation with a linear control moves  from the initial condition
 $\mathbf{{m}_0} (x)$  to the specified equilibrium
 $\mathbf{ r_1}=(-\frac{1}{\sqrt{2}},0,\frac{1}{\sqrt{2}}) .$  }
\end{figure}

Stabilization of the Landau-Lifshitz equation with affine control (\ref{eqcontrolledLL}) is illustrated in Figures~\ref{figMagnetization3d} and~\ref{figLLmcontrolTWICEICEQEQ3d} with control parameter  $k=0.5$. In Figure~\ref{figMagnetization3d}, the system dynamics are steered from $\mathbf m_0$ to 
the equilibrium point $\mathbf r_1=(-\frac{1}{\sqrt{2}},0,\frac{1}{\sqrt{2}})$.  Figure~\ref{figLLmcontrolTWICEICEQEQ3d} depicts applying the control twice in succession,   forcing the system from the equilibrium $\mathbf{r_0}$  to  $\mathbf r_2=(1,0,0)$  and then to $\mathbf r_3 = (0,0,1)$. In each case,  the state of the controlled system  converges to the specified point $\mathbf r_i $ as predicted by the analysis.

\begin{figure}[h]\hspace*{-0.3cm}
\centering\scalebox{0.544}{
    \subfloat[{\Large$\omega=1$}]{ \includegraphics{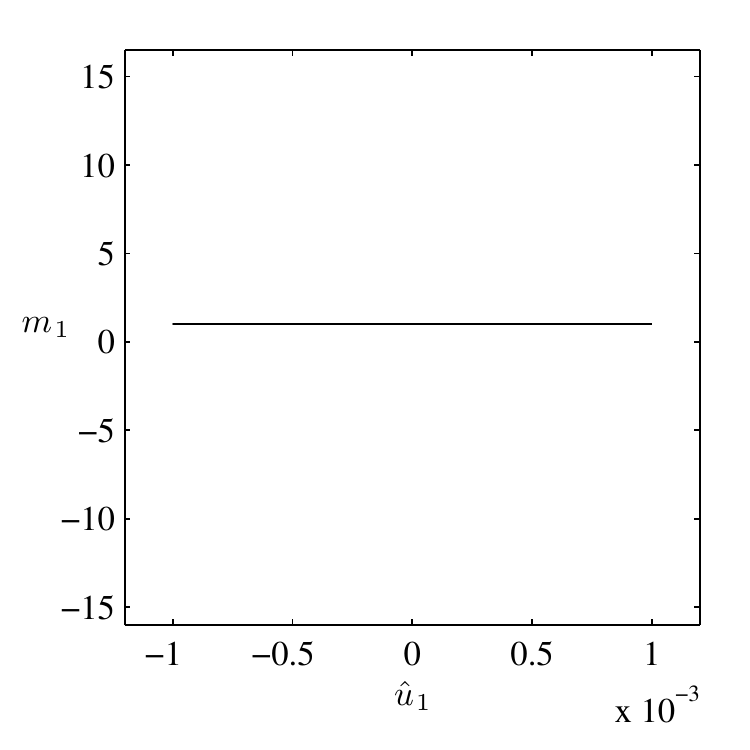}}\hspace*{0.05cm}           
  \subfloat[{\Large$\omega=0.1$}]{\includegraphics{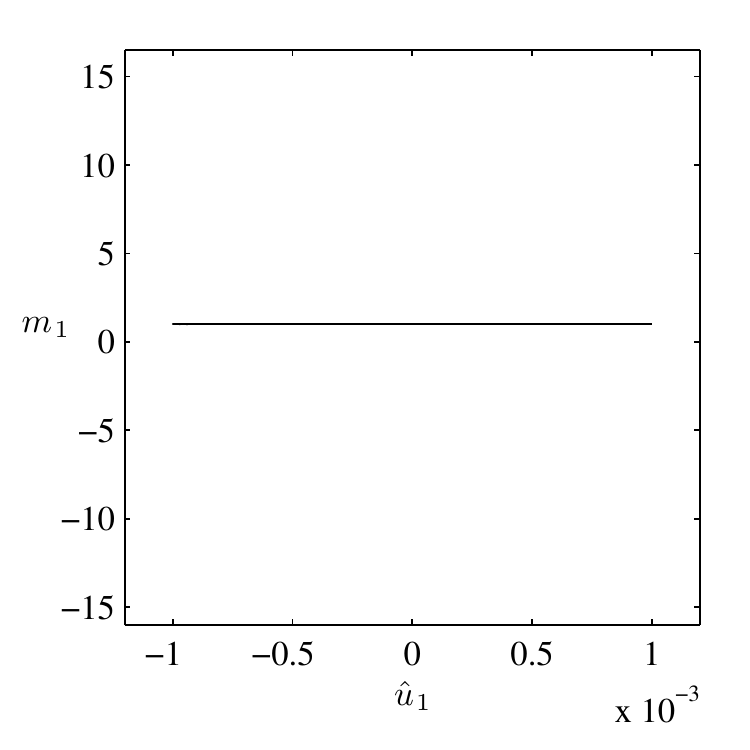}}
  \subfloat[{\Large$\omega=0.01$}]{ \includegraphics{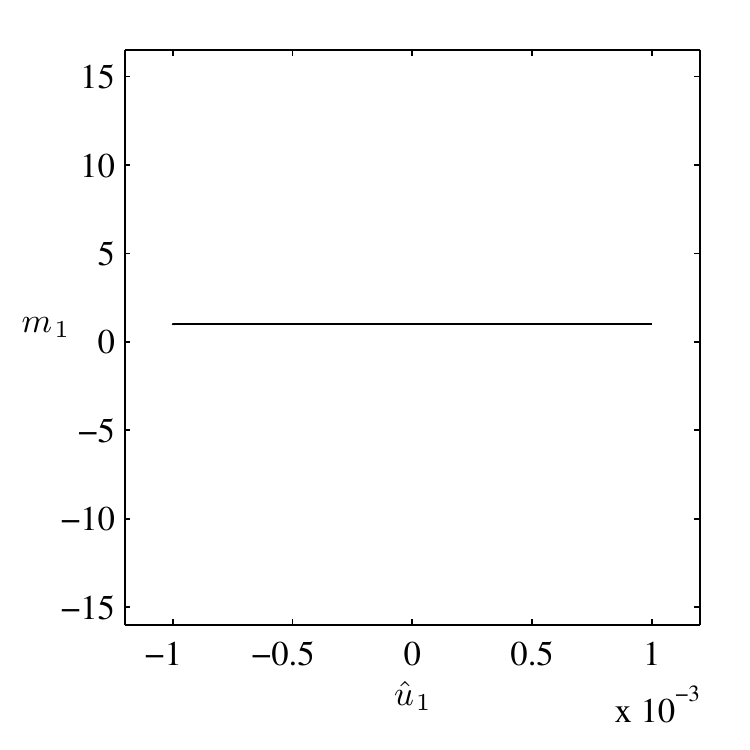}}    
    \subfloat[{\Large$\omega=0.001$}]{ \includegraphics{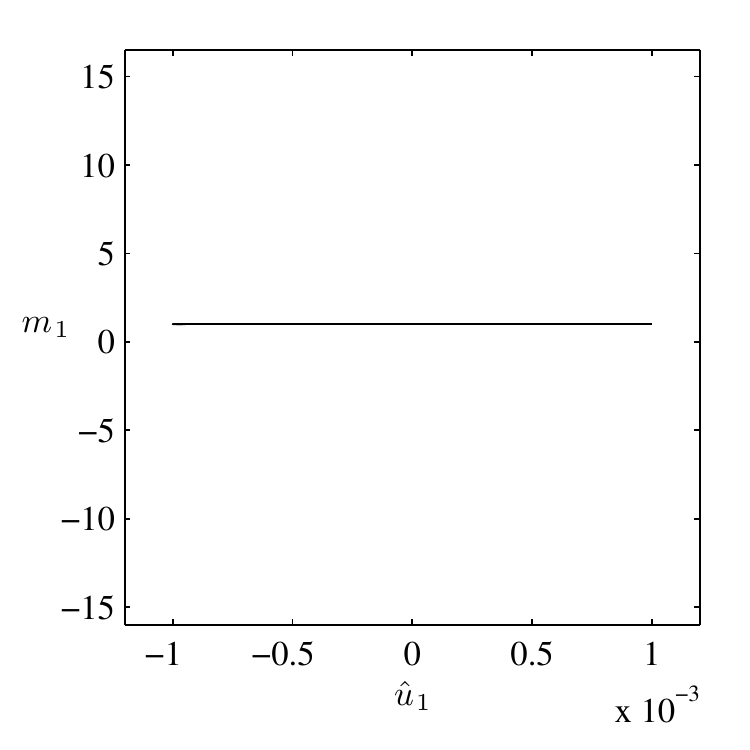}}  }  
\caption{\label{figcontrolm1LLLooping} Input--output dynamics for $m_1(x,t)$ of the controlled (nonlinear) Landau-Lifshitz equation in (\ref{eqcontrolledLL}) with $x$ fixed demonstrate the absence of persistent looping behaviour as the frequency of the periodic input $\left(0.001\cos(\omega t),0,0\right)$ approaches zero. ($L=1$, $\nu=0.02$, $\mathbf m_0(x)=\left(1,0,0\right)$, $k=0.5$, $\mathbf r=(1,0,0)$)}
\end{figure}

Adding a feedback control  so that there is only one equilibrium point also removes hysteresis from the system.  
Consider the input--output dynamics of the controlled Landau-Lifshitz equation  with periodic input $\hat{\mathbf u}(t)=\left(0.001\cos(\omega t),0,0\right)$.  The initial condition is $\mathbf m_0(x)=\left(1,0,0\right)$ and the control parameters are $k=0.5$ and $\mathbf r=(1,0,0)$. Figure~\ref{figcontrolm1LLLooping} illustrates the input--output dynamics for $m_1(x,t)$ with $x$ fixed, $L=1$ and $\nu=0.02$.  It is clear from the figure that persistent looping behaviour does not occur and hence, based on Definition~\ref{defBernstein}, the controlled Landau-Lifshitz equation in (\ref{eqcontrolledLL}) does not exhibit hysteresis. Similar behaviour is observed for $m_2(x,t)$ and $m_3(x,t)$.

\begin{figure}[h]
\centering\scalebox{0.7}{
\includegraphics{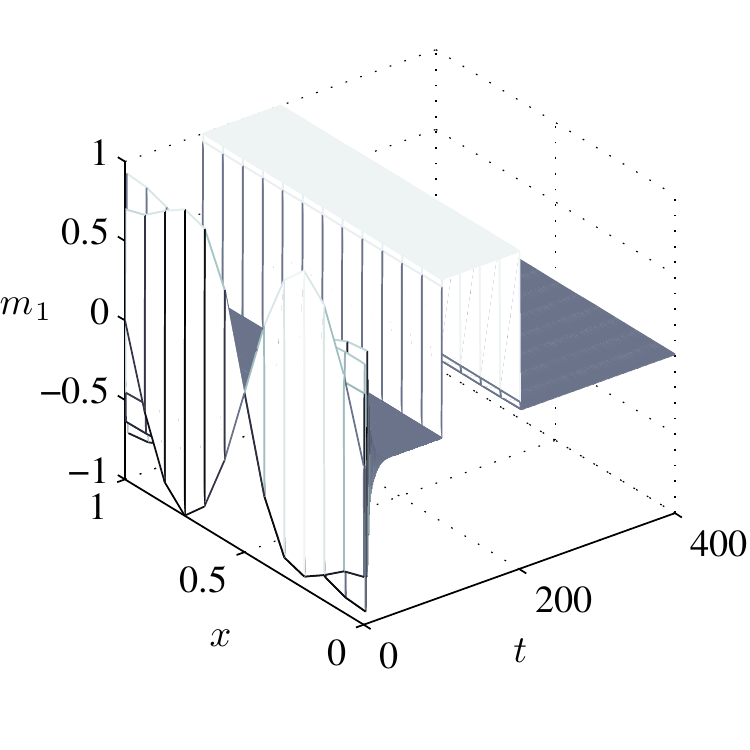}
\includegraphics{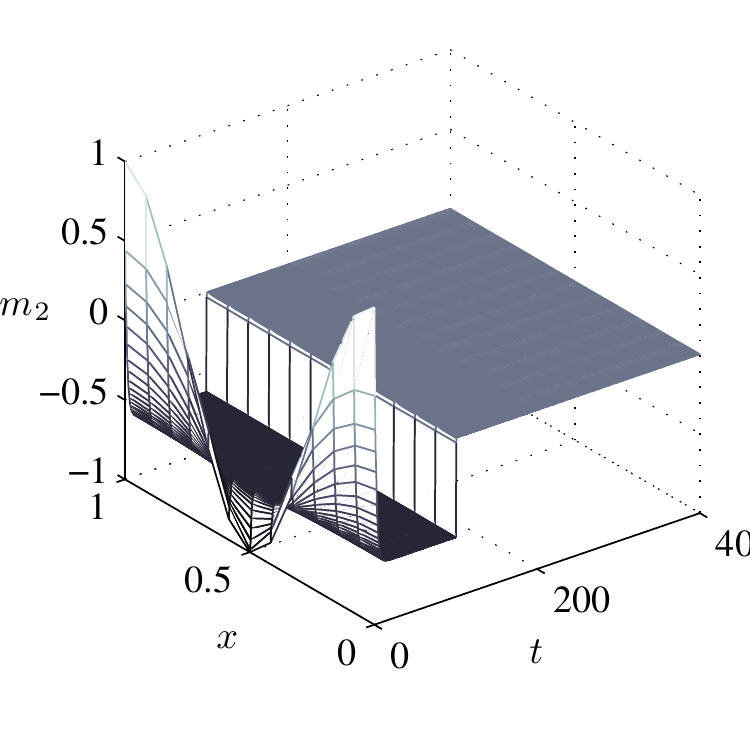}
\includegraphics{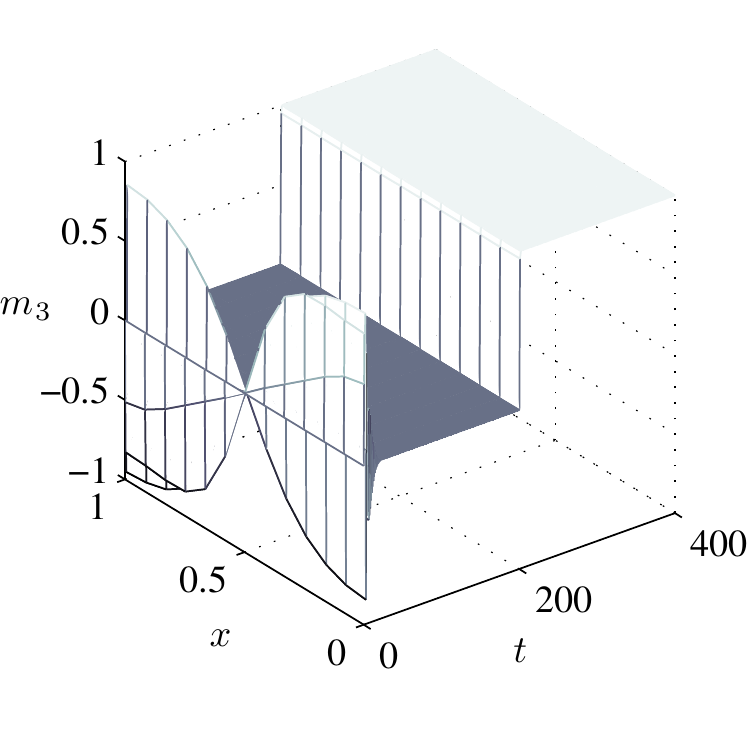}}
 \caption{\label{figLLmcontrolTwiceICEQEQ3d}  Steering magnetization between specified equilibria with a linear control. The uncontrolled magnetization moves from initial condition $\mathbf{{m}_0}$  to $\mathbf {r_0}=(0,-0.6,0)$.   Proportional control  $(k=0.5)$ with two successive values of $\mathbf r$ first forces the magnetization to $\mathbf {r_2}= (1,0,0)$ and then to $\mathbf{r_3}= (0,0,1)$. 
 }
 \end{figure}
 
\section{Conclusion}
The Landau-Lifshitz equation is a nonlinear system of partial differential equations with multiple equilibrium points. The presence of a zero eigenvalue in the linearized equation suggested a simple feedback proportional control can steer the system to an arbitrary equilibrium point. 
It was then proven that  proportional  control of the  Landau-Lifshitz equation does lead to an equilibrium point that is globally asymptotically stable in the $\mathcal L_2^3$-norm (Theorem \ref{thmr0isasymstable}) and  exponentially stable in the $H_1$-norm  (Theorem \ref{thmr0isexpstable}). 

The fact the Landau-Lifshitz equation is not quasi-linear means linearization is not guaranteed, without further analysis,
 to predict stability of the nonlinear equation \cite{alJamal2013}.  
 Moreover, since the objective of the control is to steer between equilibrium points, a linearized analysis, which only yields local results, would not predict stability of the controlled system.  Results on preservation of linearized stability require exponential stability of the linearized system; see for example  \cite[Theorem~3.3]{alJamal2013} \cite[Corollary~2.2]{Kato1995}\cite[Theorem~11.22]{Smoller1983}. 
 The fact the linearized system is exponentially stable in the $\mathcal L_2^3$-norm (Theorem~\ref{thmlinearr0isasymstable}) is encouraging, but further research is needed to determine whether the controlled nonlinear system is also exponentially stable.

\section*{Acknowledgements}

The research described in this article was supported by the Natural Sciences and Engineering Research Council of
Canada (NSERC) through grant RGPIN-6053-2015. 



\section{References}
  \bibliographystyle{elsarticle-num} 
   \bibliography{ref}

\begin{thebibliography}{10}
\expandafter\ifx\csname url\endcsname\relax
  \def\url#1{\texttt{#1}}\fi
\expandafter\ifx\csname urlprefix\endcsname\relax\def\urlprefix{URL }\fi
\expandafter\ifx\csname href\endcsname\relax
  \def\href#1#2{#2} \def\path#1{#1}\fi

\bibitem{Landau1935}
L.~Landau, E.~Lifshitz, On the theory of the dispersion of magnetic
  permeability in ferromagnetic bodies, Ukrainian Journal of Physics
  53~(Special Issue) (2008) 14--22.

\bibitem{Chow2013ACC}
A.~Chow, K.~A. Morris, Hysteresis in the {L}andau--{L}ifshitz equation, in:
  Proceedings of the American Control Conference, 2014, pp. 4747 -- 4752.

\bibitem{Morris2011}
K.~A. Morris, What is hysteresis?, Applied Mechanics Reviews 64~(5) (2011)
  050801.

\bibitem{Cowburn1999}
R.~Cowburn, D.~Koltsov, A.~Adeyeye, M.~Welland, D.~Tricker, Single-domain
  circular nanomagnets, Physical Review Letters 83~(5) (1999) 1042 -- 1045.

\bibitem{Suess2002}
D.~Suess, V.~Tsiantos, T.~Schrefl, W.~Scholz, J.~Fidler, Nucleation in
  polycrystalline thin films using a preconditioned finite element method, J.
  Appl. Phys. (USA) 91~(10) (2002) 7977 -- 9.

\bibitem{Noh2012}
S.~Noh, Y.~Miyamoto, M.~Okuda, N.~Hayashi, Y.~K. Kim, Control of magnetic
  domains in co/pd multilayered nanowires with perpendicular magnetic
  anisotropy, Journal of Nanoscience and Nanotechnology 12~(1) (2012) 428 --
  432.

\bibitem{Wiele2006}
B.~Van De~Wiele, L.~Dupra, F.~Olyslager, Memory properties in a
  {L}andau-{L}ifshitz hysteresis model for thin ferromagnetic sheets, Journal
  of Applied Physics 99~(8).

\bibitem{Yang2011}
B.~Yang, Y.~Z, Coercivity control in finite arrays of magnetic particles,
  Journal of Applied Physics 110~(10).

\bibitem{CarbouEfendiev2009}
G.~Carbou, M.~A. Efendiev, P.~Fabrie, Relaxed model for the hysteresis in
  micromagnetism, Proc. Roy. Soc. Edinburgh Sect. A 139~(4) (2009) 759--773.

\bibitem{Visintin1997}
A.~Visintin, Modified {L}andau-{L}ifshitz equation for ferromagnetism, Physica
  B: Condensed Matter 233~(4) (1997) 365 -- 369.

\bibitem{Bernstein2005}
J.~Oh, D.~Bernstein, Semilinear {D}uhem model for rate-independent and
  rate-dependent hysteresis, IEEE Trans. Autom. Control (USA) 50~(5) (2005) 631
  -- 45.

\bibitem{Guo2008}
B.~Guo, S.~Ding, {L}andau-{L}ifshitz Equations, Vol.~1 of Frontier Of Research
  with the Chinese Academy of Sciences, World Scientific, 2008.

\bibitem{Bensoussan-book}
A.~Bensoussan, Representation and Control of Infinite Dimensional Systems,
  Birkhauser, 2007.

\bibitem{Curtain1995}
R.~Curtain, H.~Zwart, An introduction to Infinite-Dimensional Linear Systems
  Theory, Vol.~21 of Texts in Applied Mathematics, Springer-Verlag, 1995.

\bibitem{LT00_1}
I.~Lasiecka, R.~Triggiani, Control Theory for Partial Differential Equations:
  Continuous and Approximation Theories, Vol.~I, Cambridge University Press,
  2000.

\bibitem{LT00_2}
I.~Lasiecka, R.~Triggiani, Control Theory for Partial Differential Equations:
  Continuous and Approximation Theories, Vol.~II, Cambridge University Press,
  2000.

\bibitem{Morris_control_handbook_rev2010}
K.~A. Morris, Control of systems governed by partial differential equations,
  in: W.~Levine (Ed.), Control Handbook, IEEE, 2010.

\bibitem{CarbouLabbe2006}
G.~Carbou, S.~Labb{\'e}, Stability for walls in ferromagnetic nanowire, in:
  Numerical mathematics and advanced applications, Springer, Berlin, 2006, pp.
  539--546.

\bibitem{Carbou2006}
G.~Carbou, S.~Labb{\'e}, Stability for static walls in ferromagnetic nanowires,
  Discrete Contin. Dyn. Syst. Ser. B 6~(2) (2006) 273--290 (electronic).

\bibitem{Carbou2011}
G.~Carbou, S.~Labbe, Stabilization of walls for nano-wires of finite length,
  ESAIM, Control Optim. Calc. Var. (France) 18~(1) (2012) 1 -- 21.

\bibitem{Jizzini2011}
R.~Jizzini, {Optimal stability criterion for a wall in a ferromagnetic wire in
  a magnetic field}, {Journal of Differential Equations} {250}~({8}) ({2011})
  {3349--3361}.

\bibitem{Labbe2012}
S.~Labbe, Y.~Privat, E.~Trelat, Stability properties of steady-states for a
  network of ferromagnetic nanowires, J. Differ. Equ. (USA) 253~(6) (2012) 1709
  -- 28.

\bibitem{Chow2015}
A.~Chow, K.~A. Morris, Control of the {L}andau-{L}ifshitz equation, Automatica
  67 (2016) 200 -- 204.

\bibitem{CoronNguyen2015}
J.-M. Coron, H.-M. Nguyen, Dissipative boundary conditions for nonlinear 1-{D}
  hyperbolic systems: sharp conditions through an approach via time-delay
  systems, SIAM J. Math. Anal. 47~(3) (2015) 2220--2240.

\bibitem{alJamal2013}
R.~al~Jamal, A.~Chow, K.~A. Morris, Linearized stability analysis of nonlinear
  partial differential equations, in: Proceedings of the 21st International
  Symposium on Mathematical Theory of Networks and Systems, 2014.

\bibitem{AM2014}
R.~al~Jamal, K.~A. Morris, Linearized stability of partial differential
  equations with application to stabilization of the {K}uramoto-{S}ivashinsky
  equation, submitted.

\bibitem{Wieser2011}
R.~Wieser, E.~Y. Vedmedenko, R.~Wiesendanger, {Indirect Control of
  Antiferromagnetic Domain Walls with Spin Current}, {Physical Review Letters}
  {106}~({6}).

\bibitem{Carbou2008}
G.~Carbou, S.~Labb{\'e}, E.~Tr{\'e}lat, Control of travelling walls in a
  ferromagnetic nanowire, Discrete and Continuous Dynamical Systems. Series S
  1~(1) (2008) 51--59.

\bibitem{Carbou2009}
G.~Carbou, S.~Labb{\'e}, E.~Tr{\'e}lat, Smooth control of nanowires by means of
  a magnetic field, Communications on Pure and Applied Analysis 8~(3) (2009)
  871--879.

\bibitem{Brown1963}
W.~Brown, Micromagnetics, no.~18 in Interscience Tracts on Physics and
  Astronomy, Wiley, 1963.

\bibitem{Gilbert2004}
T.~Gilbert, A phenomenological theory of damping in ferromagnetic materials,
  IEEE Transactions on Magnetics 40~(6) (2004) 3443 -- 3449.

\bibitem{Carbou2001}
G.~Carbou, P.~Fabrie, Regular solutions for {L}andau-{L}ifschitz equation in a
  bounded domain, Differential Integral Equations 14~(2) (2001) 213--229.

\bibitem{Alouges1992}
F.~Alouges, A.~Soyeur, On global weak solutions for {L}andau-{L}ifshitz
  equations: existence and nonuniqueness, Nonlinear Anal. 18~(11) (1992)
  1071--1084.

\bibitem{Cullity2009}
B.~Cullity, C.~Graham, Introduction to Magnetic Materials, 2nd Edition, Wiley,
  2009.

\bibitem{Lakshmanan2011}
M.~Lakshmanan, The fascinating world of the {L}andau-{L}ifshitz-{G}ilbert
  equation: an overview, Phil. Tran. Royal Society A 369 (2011) 1280--1300.

\bibitem{Banks2012}
H.~T. Banks, A functional analysis framework for modeling, estimation and
  control in science and engineering, CRC Press, Boca Raton, FL, 2012.

\bibitem{Amenda_thesis}
A.~Chow, Control of hysteresis in the {L}andau-{L}ifshitz equation, Ph.D.
  thesis, University of Waterloo (2013).

\bibitem{Luo1999}
Z.~H. Luo, B.~Z. Guo, O.~Morgul, Stability and Stabilization of Infinite
  Dimensional Systems with Applications, Communications and Control
  Engineering, Springer, 1999.

\bibitem{Kato1967}
T.~Kato, Nonlinear semigroups and evolution equations., J. Math. Soc. Japan 19
  (508--520, 1967) 508--520.

\bibitem{Michel1995}
A.~Michel, K.~Wang, Qualitative theory of dynamical systems, Vol. 186 of
  Monographs and Textbooks in Pure and Applied Mathematics, Marcel Dekker Inc.,
  New York, 1995.

\bibitem{Kato1995}
N.~Kato, A principle of linearized stability for nonlinear evolution equations,
  Trans. Amer. Math. Soc. 347~(8) (1995) 2851--2868.

\bibitem{Smoller1983}
J.~Smoller, Shock Waves and Reaction-Diffusion Equations, Vol. 258 of A Series
  of Comprehensive Studies in Mathematics, Springer-Verlag, 1983.

\end{thebibliography}






\end{document}